# Deep learning enhanced mixed integer optimization: Learning to reduce model dimensionality

Niki Triantafyllou[1,2], Maria M. Papathanasiou[1, 2*]

[1]*The Sargent Centre for Process Systems Engineering, Imperial College London, London, United Kingdom, SW7 2AZ*

[2]*Department of Chemical Engineering, Imperial College London, London, United Kingdom, SW7 2AZ*

**maria.papathanasiou11@imperial.ac.uk*

## Abstract

This work introduces a framework to address the computational complexity inherent in Mixed-Integer Programming (MIP) models by harnessing the potential of deep learning. By employing deep learning, we construct problem-specific heuristics that identify and exploit common structures across MIP instances. We train deep learning models to estimate complicating binary variables for target MIP problem instances. The resulting reduced MIP models are solved using standard off-the-shelf solvers. We present an algorithm for generating synthetic data enhancing the robustness and generalizability of our models across diverse MIP instances. We compare the effectiveness of (a) feed-forward neural networks (ANN) and (b) convolutional neural networks (CNN). To enhance the framework's performance, we employ Bayesian optimization for hyperparameter tuning, aiming to maximize the occurrence of global optimum solutions. We apply this framework to a flow-based facility location allocation MIP formulation that describes long-term investment planning and medium-term tactical scheduling in a personalized medicine supply chain.

## Introduction

Mixed Integer Programming (MIP) is a widely used modeling approach for a wide range of problems in science, engineering, and management, including planning, scheduling, control, and routing (Bertsimas & Weismantel, 2005; Floudas & Lin, 2005). In the realm of enterprise-wide optimization (EWO), many problems can be described by MIP models (Grossmann, 2005, 2012). However, a major challenge in EWO is the synchronized multi-scale decision-making across organizational functions (raw material acquisition, manufacturing, distribution, sales), geographic entities (suppliers, manufacturers, distributors, markets), and planning levels (strategic, tactical, operational) (Applequist et al., 2000; Chu and You, 2015; Dias and Ierapetritou, 2017; Laínez et al., 2012; Pistikopoulos et al., 2021; Shah, 2005; Varma et al., 2007).

EWO models can be large, multi-scale, linear, or non-linear, and hence computationally expensive to solve (Chu and You, 2015). These models, often NP-hard, are inherently combinatorial in nature due to the many discrete and logical decisions involved(Bertsimas and Weismantel, 2005; Floudas and Lin, 2005). Optimizing such problems requires effective mathematical formulations, tailored solution strategies that balance optimality with tractability, and efficient utilization of computational resources(Chu and You, 2015; Pistikopoulos et al., 2021).

To tackle these problems, different techniques have been developed to provide high-quality solutions in acceptable CPU times. These techniques involve exact algorithms (branch-and-cut, decomposition), heuristic algorithms, and matheuristics. While off-the-shelf branch-and-bound-based solvers (e.g., GUROBI, CPLEX, SCIP, Xpress, BARON, Couenne, ANTIGONE) have witnessed significant advancements, solving NP-hard MIP problems remains an open challenge. Taking advantage of the mathematical structure of large-scale models is key to being able to solve industrially-relevant MIP instances (Pistikopoulos et al., 2021). On that front, heuristic and decomposition algorithms offer a fast and effective alternative, however lacking global optimality guarantees. While the resulting solution may not always be optimal, a well-designed heuristic method often yields a good candidate solution. However, heuristic and decomposition algorithms often rely on problem-specific knowledge, necessitating extensive trial and error by the modeler (Khalil et al., 2017a). These approaches include relaxed and aggregated model formulations, hierarchical methods where information flows through subproblems, iterative methods through the use of integer cuts until the integrality gap is closed between subproblems, and full space decomposition algorithms (Maravelias and Sung, 2009).

Full space decomposition methods can be categorized into algorithms targeting complicating constraints and those addressing complicating variables (Allman et al., 2019; Conejo et al., 2006; Daoutidis et al., 2019). Decomposition algorithms addressing complicating constraints include the Dantzig-Wolfe decomposition for linear programming (LP) (Allen et al., 2023; Sun et al., 2018; Wu et al., 2020), its MIP extension called Branch-and-price (Aguayo et al., 2019; Allman and Zhang, 2021; Cóccola et al., 2015; Gunpinar and Centeno, 2016), Lagrangean decomposition (Oliveira et al., 2013; Tautenhain et al., 2021; Terrazas-Moreno et al., 2011; van den Heever et al., 2001; Yang et al., 2020; You et al., 2011), augmented Lagrangean decomposition (Gharaei et al., 2019; Nishi et al., 2008), and the Alternating Method of Multipliers (ADMM) (Chang et al., 2015; Garcia-Herreros et al., 2016; Labbé et al., 2019). On the other hand, algorithms handling complicating variables include Benders decomposition (R. Cory Allen et al., 2023; Bhosekar et al., 2021; Fontaine and Minner, 2014; Kim and Zavala, 2018; Lara et al., 2018; Leenders et al., 2023; Saharidis and Ierapetritou, 2009; You and Grossmann, 2013), generalized

Benders decomposition (Gharaei et al., 2020; Mitrai and Daoutidis, 2022) and bilevel decomposition (Calfa et al., 2013; Elsido et al., 2019; Terrazas-Moreno and Grossmann, 2011; You et al., 2011).

### 1.1 Mixed integer programming background

Many problems in EWO can be described by MIP models (Grossmann, 2012). In this study MIP and mixed integer linear programming (MILP) are used interchangeably. MIP models seek to find an optimal solution $x^*$ in a discrete set $\mathcal{F}$ that minimises or maximises an objective function $f(x)$ for all $x \epsilon \mathcal{F}$. Given matrices $\boldsymbol{A} \in \mathbb{Z}^{n \times m}$ , $\boldsymbol{B} \in \mathbb{Z}^{n \times k}$ and vectors $\boldsymbol{b} \in \mathbb{Z}^{n}$ , $\boldsymbol{c} \in \mathbb{Z}^{m}$, $\boldsymbol{d} \in \mathbb{Z}^{k}$, a mixed integer programming (MIP) problem can be formulated as (Bertsimas and Weismantel, 2005):

$$
\begin{aligned}
\min \quad & \boldsymbol{c}^T \boldsymbol{y} + \boldsymbol{d}^T \boldsymbol{x} \\
\text{s.t} \quad & \boldsymbol{A}\boldsymbol{y} + \boldsymbol{B}\boldsymbol{x} \leq \boldsymbol{b} \\
& \boldsymbol{y} \in \mathbb{Z}_+^m, \quad \boldsymbol{x} \in \mathbb{R}_+^k
\end{aligned}
$$

### 1.2 Machine Learning for Mixed Integer Programming

Machine learning (ML) can be leveraged to tackle MIP problems by employing various techniques and algorithms. State-of-the-art algorithms often rely on manually crafted heuristics to make decisions, and to perform expensive computations (Bengio et al., 2021). In this context, ML is a promising alternative to enhance these decisions in a more efficient manner. Machine learning can contribute to MIP improvement in two primary ways: (a) by replacing expensive computations with fast approximations, and (b) by addressing algorithmic decisions that heavily depend on expert knowledge and intuition, which may result in suboptimal outcomes. Two excessive reviews on machine learning for discrete optimization are given by Bengio et al. (2021) and Zhang et al. (2023).

There are several approaches that machine learning can be used for combinatorial optimization. Machine learning models trained with precomputed MIP instances can be used to directly replace MIP problems (Abbasi et al., 2020; Larsen et al., 2018; Vinyals et al., 2015). Alternatively, machine learning can be used to augment exact methods with additional problem-specific information (Ammari et al., 2023; Bertsimas and Margaritis, 2023; Bertsimas and Stellato, 2022; Beykal et al., 2022; Chu and You, 2014; Goettsch et al., 2020; Kilwein et al., 2023; Xavier et al., 2020). Machine learning can also be used to enhance decomposition algorithms (Abbas and Swoboda, 2022; R. Cory Allen et al., 2023; Mitrai and Daoutidis, 2024a) or to decide whether to decompose a problem or not (Basso et al., 2020; Kruber et al., 2017; Mitrai and Daoutidis, 2024b). These approaches acknowledge that the combination of machine learning with optimization techniques can lead to more effective solutions. Finally, machine learning models can also be used alongside optimization algorithms, where they are iteratively called through the optimization algorithm with the scope of improved algorithmic performance (Bengio et al., 2021). This approach has been used to enhance the branch-and-cut algorithm by learning for branching variable selection (Gasse et al., 2019; Hottung et al., 2019; Khalil et al., 2016; Lodi and Zarpellon, 2017; Zarpellon et al., 2021), node selection (Yilmaz and Yorke-Smith, 2021), primal heuristics (Khalil et al., 2017b; Nair et al., 2020), and cutting plane selection (Baltean-Lugojan et al., 2018; Huang et al., 2022; Li et al., 2023; Paulus et al., 2022).

Solving MIP models, efficiently, remains a challenge, particularly as enterprise-wide optimization problems become more complex in the global setting. In this context, machine learning offers implicit knowledge that can complement the explicit knowledge gained from mixed integer programming. ML offers a powerful tool for automatically developing effective heuristics from a dataset of mixed integer programming (MIP) instances. This can be particularly valuable in applications where numerous instances of the same problem must be solved and while the constraints and variables of each instance change significantly, but they do so in a predictable manner. In real-world scenarios, problem instances often vary and belong to different distributions. There's a need for solvers that can adapt specifically to these unique problem-specific distributions, allowing for high quality solutions in fast computational times.

In this work, we present a deep-learning enhanced approach to decompose large-scale MILP problems by approximating complicating variables and identifying the problem's active dimensions. We show that deep learning constructed heuristics substantially improve the performance of standard off-the-shelf MIP solvers. Our approach enhances and automates mixed integer programming with expert intuition retrieved from the specific structure of the problem at hand.

# 1. Methodology

We present a novel approach to address computational complexity in Mixed-Integer Programming (MIP) models. The objective is to develop and use a deep learning model based on neural networks that approximates complicating binary variable(s) $y_m \in \{0,1\},\ \forall\, m \in \boldsymbol{M}$ of the original MIP problem and therefore reducing its dimensionality $y_{m^*} \in \{0,1\},\ \forall\, m^* \in \boldsymbol{M}^* \subseteq \boldsymbol{M}$. The framework starts by taking the MIP model of interest and generating a set of random *N* problem instances drawn from different distributions. This is done to ensure better generalization. Following this, the MIP is solved for each instance with a standard off-the-shelf solver, creating a comprehensive input-output dataset, which is used for training the deep learning model by learning the underlying problem-specific structure. To enhance the framework's performance, we employ Bayesian optimization for hyperparameter tuning, aiming to maximize the occurrence of global optimum solutions.

Once trained, the deep learning classifier is employed to make predictions for a target problem instance, approximating the complicating binary variable(s). The resulting approximation can be leveraged (a) *indirectly*, by identifying the active dimensions associated with the complicating variables for each problem instance; and (b) *directly*, by fixing the complicating variables. The indirect approach, demonstrated through a case study, leads to a reduced MIP model with significantly fewer constraints and binary variables than the original model. This approach is particularly relevant for complex problems such as facility location, capital budgeting, and variations thereof, including supply chain optimization, investment planning, energy production, and distribution, as well as project portfolio selection. On the other hand, the direct fixing method decomposes the problem, allowing its solution in blocks, suitable for problem classes such as vehicle routing, inventory management, and process scheduling, where all dimensions are active. In both cases, the resulting reduced MIP is optimized with a commercial branch-and-cut solver, and the optimal solution is obtained. The aim of our methodology is to enhance the branch-and-cut algorithm with problem-specific information, without relying on expert knowledge.

The presented methodology (Figure 1) comprises three main components: (a) synthetic data generation, (b) deep-learning classifier training and hyper-parameter tuning with Bayesian optimization for automatic construction of heuristics, and (c) combination of the deep learning heuristics with the branch-and-cut.

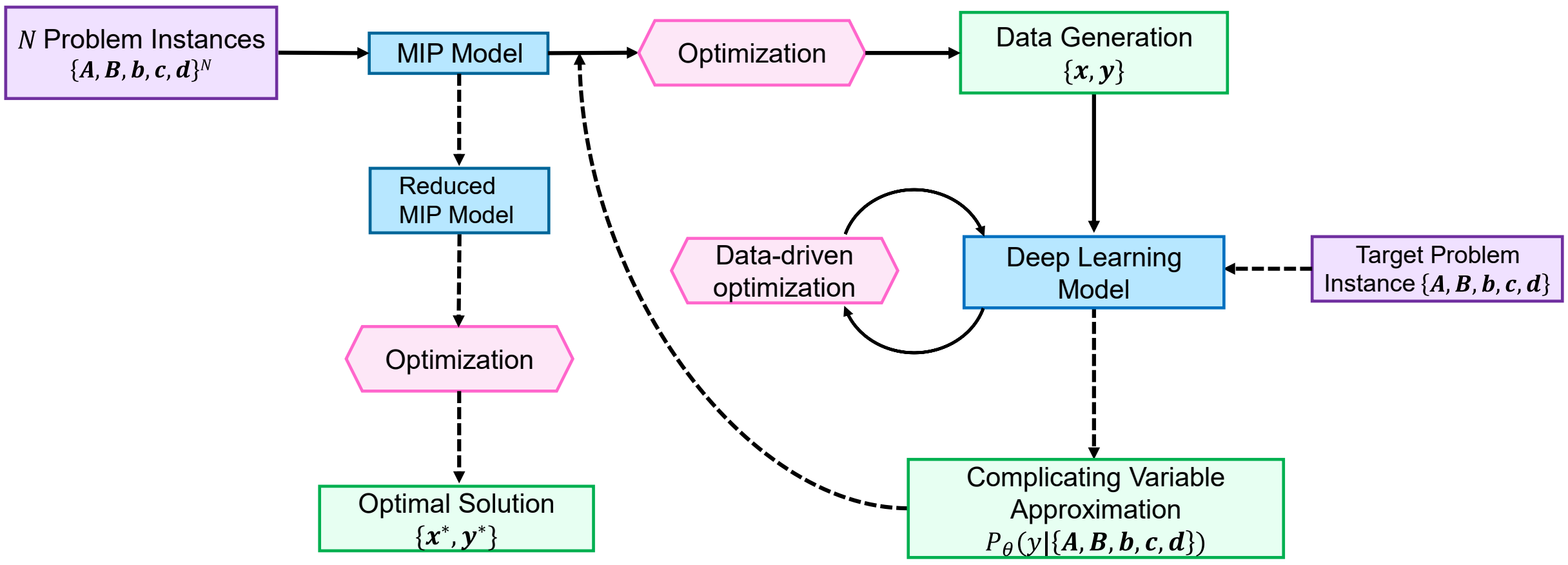

Figure 1. Proposed decomposition framework for Mixed Integer Programming (MIP) models. Solid lines (–) signify the steps undertaken once during the deep learning model training, while dashed lines (--) represent steps that can be repeated as needed to solve any instance of the MIP post-model training.

## 2.1 Synthetic Data Generation

The first step focuses on the generation of a comprehensive dataset used to train the deep learning classifier. This dataset is created through the solution of a representative number of instances of the original MIP model with an off-the-shelf MIP solver. The parameter instances of the MIP involve the $\{A, B, b, c, d\}$ of the MIP representation in section 1.1. The number of instances needed to train a robust classifier depends on the problem at hand. The incorporation of synthetic data generation is imperative, particularly in cases where real-world data is scarce or entirely absent. In cases where real-world data is available, it can either be used independently or augmented with synthetic data to enhance the uniformity of the dataset(de Melo et al., 2022).

The objective of the synthetic data generation step is to establish a direct link between uncertain parameters and the complicating variables of the MIP. Algorithm 1 outlines the proposed step-by-step procedure for synthetic data generation. This approach ensures that any uncertainty related to the variability of the critical parameters is captured and incorporated into the training dataset.

**Algorithm 1: Synthetic Data Generation**

**Inputs:**

- Mixed Integer Optimization model $MIP$
- A target binary complicating variable $y_m \in \{0,1\}, \; \forall \, m \in M$
- Subset $n \subseteq N_{par}$ of uncertain parameters $\boldsymbol{P} = \{p_1, p_2, \ldots, p_n\}$ of $\{\boldsymbol{A}, \boldsymbol{B}, \boldsymbol{b}, \boldsymbol{c}, \boldsymbol{d}\}$ and their upper and lower bounds $[\boldsymbol{P}_{min}, \boldsymbol{P}_{max}]$
- Number of instances $N_{ins}$ for the critical uncertain parameters
- Number of different probability distributions $N_{dist}$ for each instance
- Number of problem instances per probability distributions $N_{ins,dist}$;

Initialize the feature vector $InstanceSet = \emptyset$<br>
$i = 1$<br>
**while** $i \leq N_{ins}$ **do**<br>
    **for** $j = 1$ to $N_{dist}$ **do**<br>
        **for** $k = 1$ to $N_{ins,dist}$ **do**<br>
            Generate randomized problem instance $\{\boldsymbol{A}, \boldsymbol{B}, \boldsymbol{b}, \boldsymbol{c}, \boldsymbol{d}\}_{i,j,k} \leftarrow Prob_j(\boldsymbol{P}_{min}, \boldsymbol{P}_{max})$ ;<br>
            Append $InstanceSet(\{\boldsymbol{A}, \boldsymbol{B}, \boldsymbol{b}, \boldsymbol{c}, \boldsymbol{d}\}_{i,j,k})$ ;<br>
        **end for**<br>
    **end for**<br>
    $i \leftarrow i + 1$<br>
**end while**<br>
Obtain $InstanceSet$<br>
$\boldsymbol{DataSet} = \emptyset$<br>
**for** $i = 1$ to $|InstanceSet|$ **do**<br>
    Solve $MIP(InstanceSet_i)$ using a commercial solver;<br>
    **if** $MIP(InstanceSet_i)$ $is\ feasible$ **then**<br>
        Obtain value of complicating variable $y_{m,i} \; \forall \, m \in M$ ;<br>
        Append $DataSet(InstanceSet_i, \; y_{m,i} \in \{0,1\}^m; \; infeasible = 0)$ ;<br>
    **else**<br>
        Append $DataSet(InstanceSet_i, \; y_{m,i} = 0 \; \forall \, m; \; infeasible = 1)$ ;<br>
    **end if**<br>
**end for**

**Output**: $\boldsymbol{DataSet} = \{(InstanceSet, \; y_m\}_{i=1}^{N}$ where $y_m \in \{0,1\}^{M+1}$

### 2.1.1 Complicating variable identification

To create the dataset, firstly, the complicating variables of the optimization problem at hand need to be identified. In the case of mixed integer linear programming (MILP) problems, the integer variables are usually considered to be the complicating variables (Conejo et al., 2006). This is because the remaining continuous problem can be decomposable by blocks. For example, in the context of long-term investment planning problems, integer variables often govern investment decisions, while continuous variables pertain to operational decisions.

### 2.1.2 Problem instance generation

Upon identifying the set of target *complicating variables*, the critical *uncertain parameters* need to be selected. These can be multi-dimensional and are typically bounded within the ranges of interest. Different realizations of those parameters lead to different outcomes for the complicating variables. It should be noted that the uncertain parameters are the features of the problem. Typically, uncertain parameters are multidimensional, hence they are aggregated in a suitable manner to represent a specific number of features for each instance. For example, if the uncertain parameter is a multi-period (discretized daily) demand profile, each time point (day) within the profile or sums of time points (months) can be formulated as distinct features. To ensure a representative dataset, it is important to choose an appropriate number of problem instances $N_{ins}$ for randomized realizations of the uncertain parameters.

Ensuring comparable performance between the instances used for learning (training) and those unseen (testing and validation) is commonly referred to as generalization in machine learning. Usually, deep learning algorithms generalize well to examples within the same distribution but often encounter challenges to generalize to out-of-distribution instances (Bengio et al., 2021). To address this and guarantee dataset impartiality and classifier robustness to parameter uncertainty, we introduce an additional layer of instance generation. A number of probability distributions $N_{dist}$ is selected to sample the parameter realizations. The number and the type of probability distributions are defined by the modeler, and it is dependent on the optimization problem and any prior knowledge regarding the uncertain parameters. Finally, we define the number of instances to be generated for each probability distribution $N_{ins,dist}$. Consequently, the total number of randomly generated instances is calculated as $N_{ins} \times N_{dist} \times N_{ins,dist}$. During the generation of each instance, the features are systematically stored in the $InstanceSet$. This approach ensures the generation of a diverse dataset capable of training a classifier that incorporates uncertainties introduced through random sampling.

### 2.1.3 Imitation learning

Imitation learning is a class of algorithms that seek to learn the policy of an expert (oracle) through examples of expert behaviour. In our approach, this is formulated as a supervised learning problem. In supervised learning, the goal is to teach a model to predict labels for input instances based on training data. The optimal labels are usually learned through an expert. In this context, the expert that provides the optimal label is the solver that solves the mixed integer optimization model to global optimality. To ensure the solver finds the best possible labels, the optimality gap is set to 0.01%. The optimality gap is a measure of how close the relaxed solution of the optimization problem is to the solution of the original optimization problem. By setting the optimality gap to 0.01%, the goal is to minimize the Bayes error by making sure that the oracle always provides the best solution possible. The Bayes error is the lowest possible error that any classifier can achieve for a given distribution of data. To generate the label values, the MIP model is solved for all instances in the $InstanceSet$ using commercially available solvers.

The binary labels of the $DataSet$ are the dimensions $M = \{m_1, m_2, m_3 \dots, m_n\}$ of the complicating variable $y_m \in \{0,1\}, \ \forall\, m \in M$. If the binary variable $y_{m_n}$ is active for a certain dimension $m_n$, which means that the corresponding label $m_n$ is 1, otherwise it is 0 (Algorithm 1). An additional binary label is added to note instance infeasibility (*infeasible*). If an MIP instance is infeasible, then the infeasibility label is noted as 1, otherwise is 0. This can be critical, especially in the case of large problem instances, where the solver may require significant computational time (several hours or even days) to determine whether a problem is infeasible. The dataset is split into 80% training, 10% testing, and 10% validation sets.

## 2.2 Deep learning

Here the most complex step of the optimization process is replaced by a predictive multi-label classification model, in which each distinctive dimension of the complicating variable represents a class label. In this work, we harness the predictive capabilities and the ability to extract complex patterns from multiple dimensions (features) of deep learning to build the multi-label classification model. The inputs to the classification model are the uncertain parameters of the MIP model (aka MIP problem instances) $\{\boldsymbol{A}, \boldsymbol{B}, \boldsymbol{b}, \boldsymbol{c}, \boldsymbol{d}\}_i \in InstanceSet$ for $i = 1,2, \dots N$. The labels are the binary complicated variables $y_{m,i} \in \{0, 1\}^{|M|+1}$, where the total number of labels is the cardinality of the dimension set $M$ plus one extra label for infeasibility, therefore $|M| + 1$.

### 2.2.1 Multi-label classification

Multi-label classification is used to predict the active dimensions of the complicating variable. This allows each input problem instance to be linked with multiple non-exclusive class labels (Tsoumakas and Katakis, 2007). This becomes particularly crucial when dealing with multi-dimensional complicating binary variables (e.g. candidate manufacturing facilities in planning problems or candidate products in scheduling problems), as they often have more than one active dimension for each problem instance. In this case, we employ a specific case of multi-label classification, whereby the infeasibility class is exclusive. This means that when the *infeasible* class is active, it dictates the inactivity of all other classes (dimensions).

### 2.2.2 Neural network architecture and training

Neural networks can be configured to map the input instances of the MIP problem ($\boldsymbol{INS}$) to the complicating variable $\boldsymbol{y_m}$ by approximating a function $\boldsymbol{y_m} = f_m(\{\boldsymbol{A}, \boldsymbol{B}, \boldsymbol{b}, \boldsymbol{c}, \boldsymbol{d}\}; \boldsymbol{w_m})$ for every label $m \in$ $|M| + 1$. Here we investigate and compare the performance of two different architectures: (a) feedforward neural networks (also known as multilayer perceptron or ANN), and (b) convolutional neural networks (ConvNet or CNN).

**(a) Feed-forward neural networks**

Input data undergoes transformations through multiple layers (Goodfellow et al., 2016). Each layer applies an affine transformation to the input vector (($\{\boldsymbol{A}, \boldsymbol{B}, \boldsymbol{b}, \boldsymbol{c}, \boldsymbol{d}\}$), followed by the non-linear activation function applied element-wise. The rectified linear unit (ReLU) is used as the activation function (LeCun et al., 2015), while the output of one layer becomes the input for the next layer. If the output of layer $l$ is denotes as $\boldsymbol{a}^{(l)}$ (before activation) and $\boldsymbol{o}^{(l)}$ (after activation), then the transformations from layer $l - 1$ to layer $l$ is given by Eqs. (1), (2):

$$\boldsymbol{a}^{(l)} = \boldsymbol{W}^{(l)}\boldsymbol{o}^{(l-1)} + \boldsymbol{b}^{(l)} \tag{1}$$

$$\boldsymbol{o}^{(l)} = \mathrm{ReLU}(\boldsymbol{a}^{(l)}) \tag{2}$$

Where, $\boldsymbol{W}^{(l)}$ is the weight matrix for layer $l$ and $\boldsymbol{b}^{(l)}$ is the bias vector for layer $l$.

The output of the last layer, denoted as $\boldsymbol{a}^{(output)}$, is used for predictions. The sigmoid activation function is then applied element-wise to $\boldsymbol{a}^{(output)}$ to obtain the final predictions Eq. (3):

$$\hat{\boldsymbol{y}} = \mathrm{sigmoid}(\boldsymbol{a}^{(output)}) \tag{3}$$

In multi-label classification, the number of output layers is dictated by the number of target labels $|M| + 1$. The sigmoid activation function is used in the output layers to transform the output values into a range between 0 and 1, representing probabilities. The sigmoid function allows each output layer to make independent predictions, indicating the likelihood of a specific class (i.e. dimension of the complicating variable) being present (Goodfellow et al., 2016). This independence is crucial in the case where multiple dimensions can be active for a single MIP instance.

The weights ($\boldsymbol{W}$) and biases ($\boldsymbol{b}$) in each layer are learned (optimized) through stochastic gradient descent, which minimizes a chosen loss function $L_{NN}$ (LeCun et al., 2015), using the backpropagation algorithm to compute gradients efficiently. Each label is treated as an independent binary classification task and binary cross entropy loss (Eq (4)) is chosen as loss function to measure the dissimilarity between true ($y_{m,i}$) and predicted binary outcomes ($\hat{y}_{m,i}$).

$$L_m = -\frac{1}{n}\sum_{i=1}^{n}(y_{m,i} \cdot \log \hat{y}_{m,i} + (1 - y_{m,i}) \cdot \log(1 - \hat{y}_{m,i})) \tag{4}$$

The overall loss for the neural network is the sum of losses across all labels given by Eq (5):

$$L_{NN} = \sum_{m \in |M|+1} L_m \tag{5}$$

**(b) Convolutional neural networks**

Convolutional Neural Networks (CNNs) specialize in processing data with a grid-like topology, capturing spatial and temporal relationships in input data(Goodfellow et al., 2016). A CNN typically consists of convolutional, pooling, and fully connected layers. The typical CNN structure is organized into stages, with the initial stages comprising convolutional and pooling layers. The former apply convolutions that are linear operations that involve filters (kernels) sliding across the input data, capturing patterns and features(LeCun et al., 2015), arranged in feature maps. These features maps share weights from layer to layer, promoting the detection of patterns with shared weights. The convolutional operation is followed by a ReLU activation function to capture non-linearities. Pooling layers play a crucial role in merging semantically similar features by reducing the dimension of the feature maps and introducing invariance to small shifts and distortions(LeCun et al., 2015).

Let's denote the output of a convolutional layer as $\boldsymbol{C}^{(l)}$. The convolutional operation is followed by the ReLU activation function (Eq (6)):

$$\boldsymbol{C}^{(l)} = \mathrm{ReLU}(\mathrm{Conv}(\boldsymbol{C}^{(l-1)}, \boldsymbol{F}^{(l)}) + \boldsymbol{b}^{(l)}) \tag{6}$$

Where $\mathrm{Conv}(\cdot)$ is the convolution operation, $\boldsymbol{F}^{(l)}$ is the set of filters for layer $l$ and $\boldsymbol{b}^{(l)}$ is the bias term for layer $l$

Next, the max pooling operation is performed by computing the maximum output from specific neighbourhoods of the feature map based on the kernel and stride size. The pooling layer is denoted as $\boldsymbol{P}^{(l)}$ given by Eq (7):

$$\boldsymbol{P}^{(l)} = \mathrm{MaxPooling}(\boldsymbol{C}^{(l)}) \tag{7}$$

After the convolutional layers, the output feature map is flattened into a vector (Eq (8)):

$$\boldsymbol{V} = \mathrm{Flatten}(\boldsymbol{C}^{(last)}) \tag{8}$$

Two or three stages of convolution, non-linearity, and pooling are stacked to capture hierarchical representations of the input data. These are followed by fully-connected layers and the sigmoid output layer as described in the feed-forward neural network section. Additionally, dropout is used in the fully connected layers as a regularization technique to reduce overfitting.

### 2.2.3 Data imbalance

Class imbalance can often challenge multi-label classification (Pant et al., 2019). Multi-label datasets consist of a variety of classes, which may result in an uneven distribution of these classes throughout the dataset. This can lead to: (i) inter-class imbalance, where some labels are much more frequently used than others, and (ii) inner-class imbalance, where there is an imbalance between positive and negative examples for each label (Liu et al., 2021; Pirizadeh et al., 2023). Inner-class imbalance is usually inversely proportional to the frequency of each label; thus it is more frequently encountered where there exist rare labels. In multi-label datasets, negative classes tend to be over-represented, while positive classes are often under-represented. This imbalance poses a significant challenge when training classification models, as it reduces the ability to effectively learn from and predict minority classes. In the context of building a classifier to identify the active dimensions of a complicating variable in a MIP problem, the class/dimension imbalance is an inherent challenge of the problem. Some dimensions are more likely to be active than others or some dimensions are active only if another dimension is active. Addressing class imbalance in multi-label classification is essential for improving the model's performance. Some common strategies include (a) undersampling the majority class or oversampling the minority class (e.g. SMOTE)(Liu et al., 2021), (b) synthetic data generation for minority classes, (c) cost-sensitive learning that takes into consideration misclassification costs (Ling & Sheng, 2010)(Ling & Sheng, 2010)(Ling & Sheng, 2010), (d) ensemble methods, such as bagging or boosting (Khoshgoftaar et al., 2011), and (e) active learning (Aggarwal et al., 2021; Ertekin et al., 2007).

#### 2.2.4 Evaluation metrics

A valuable metric for assessing the performance of a classifier is the confusion matrix. However, in the context of multi-label classification, the traditional confusion matrix used for multi-class classification is undefined. To address this, here we use the multi-label confusion matrix with one extra column for the non-predicted labels (NPL). NPL accounts for labels that the classifier failed to predict, while one extra row for the no true labels (NTL) is also incorporated (Heydarian et al., 2022). To ensure that the performance of the classifier is accurately assessed and to account for cases where the classifier predicts a subset of the correct values but not all, both partial and total correctness metrics are used. Total correctness is assessed through sample-level accuracy, while partial correctness is assessed using metrics such as hamming loss, Jaccard index, precision, recall, and F1-score. For the definitions and further details on the evaluation metrics used, refer to Appendix B.

#### 2.2.5 Hyperparameter tuning

Deep neural networks have hyperparameters that are used to control their training. Hyperparameters are not learned and optimized like the internal parameters (weights and biases) of the neural networks, but they are manually set before training a model (Feurer and Hutter, 2019). Examples of hyperparameters include the learning rate, number of hidden layers, number of neurons, batch size, number of training epochs, dropout rate, etc. Effective hyperparameter tuning is crucial for optimal model performance. Hyper-parameter tuning techniques such as random search, grid search, and Bayesian optimization aim to discover improved hyperparameter combinations (Feurer and Hutter, 2019).

In our approach, we employ Bayesian optimization based on Gaussian process ($GP$) surrogates for hyperparameter tuning as seen in Algorithm 2 (Wu et al., 2019). An initial $GP$ model is trained based on 5 samples generated through the Sobol quasi-random sequence. Then the $GP$ model is iteratively updated through the upper confidence bound (UCB) acquisition function, balancing exploration and exploitation(Brochu et al., 2010; Rasmussen and Williams, 2005).

Viewing the neural network as a black box function, Bayesian optimization maximizes a chosen evaluation metric—in this case, sample-level accuracy. Sample-level accuracy accounts for total correctness, which means that it rewards only the set of dimensions (labels) of the complicating variable that leads to the global optimum solution. By maximizing the sample-level accuracy, we maximize global optimality and penalize local optimum and infeasible combinations of the predicted complicating variable dimensions. With this approach, the neural network is trained with the optimal hyperparameters to favor global optimal solutions.

**Algorithm 2: Bayesian Optimization for Global Optimality Maximization**

**Input:** $DataSet$ from Algorithm 1, $NN(\cdot)$ type, hyper-parameters $\boldsymbol{\theta} = \{\theta_1, \theta_2, \dots, \theta_n\}$, sample-level accuracy function $A(\boldsymbol{\theta})$, maximum number of evaluations $maxiter$;

1 Sample 5 hyperparameter combinations $\boldsymbol{\theta}_0$ based on the Sobol quasi-random sequence;

2 Train $NN(\boldsymbol{\theta}_0, DataSet)$ given the $DataSet$ and $\boldsymbol{\theta}_0$ and obtain the accuracy $A(\boldsymbol{\theta}_0)$;

3 **BO initialization:** train $GP^0$ with an initial data set $D_0 = (\boldsymbol{\theta}_0, A(\boldsymbol{\theta}_0))$;

4 **for** $i = 1$ to $maxiter$ **do**

5 Find $\boldsymbol{\theta}_i$ by optimizing the upper confidence bound acquisition function over the $GP^i$, where $\boldsymbol{\theta}_i \leftarrow argmin_{\boldsymbol{\theta}}\ UCB(\boldsymbol{\theta}|D_{1:i-1})$;

6 Train the $NN(\boldsymbol{\theta}_i, DataSet)$ and obtain the accuracy $A(\boldsymbol{\theta}_i)$;

7 Augment the data set $D_{1:i} = \{D_{1:i-1}, (\boldsymbol{\theta}_i, A(\boldsymbol{\theta}_i))\}$ and update the $GP^i$;

8 **end for**

**Output**: Trained neural network $NN(\boldsymbol{\theta}^*, DataSet)$

### 2.2.6 Reduced mixed integer optimization model

Once the trained neural network is obtained, it can be used to configure the complicating variable in the MIP model by reducing its dimensionality indefinitely. The decomposition procedure is given in Algorithm 3. Given a target MIP parameter instance, the trained neural network approximates the complicating variable as a probability for each class $y_m = NN(\{\boldsymbol{A}, \boldsymbol{B}, \boldsymbol{b}, \boldsymbol{c}, \boldsymbol{d}\}^*; k_{prob}) \ \forall\, m \in \mathbf{M}$ and therefore the active dimensions of $\boldsymbol{M}$. To decide which dimensions are active based on these probabilities, a threshold $k_{prob}$ (commonly set at 0.5) is used. Lowering $k_{prob}$ makes the algorithm more conservative, increasing computational time while reducing the likelihood of infeasible solutions.

There are two approaches that can be followed, (a) either the predicted complicating variable is fixed in the $MIP$ or (b) the active dimensions of the complicating variable $m \in \mathbf{M}$ are used to reduce the dimensionality of set $\mathbf{M}$ to $\boldsymbol{M}^* \subseteq \boldsymbol{M}$. Both approaches lead to a reduced $MIP$ model, called $RMIP$. The reduced $MIP$ ($RMIP$) is now solved for the target instance using a commercial solver, and finally, the optimal solution is obtained. As this is a hierarchical decomposition approach, there are no mathematical guarantees of optimality.

Training a perfect classifier is a challenge, hence the trained neural network may sometimes lead to suboptimal or infeasible solutions. To address this limitation, we adopt a probabilistic approach by obtaining the predicted probabilities for all the labels. Subsequently, we introduce a user-specified threshold for these probabilities, allowing us to identify all labels (dimensions) that surpass this threshold. This user-defined threshold serves as an important parameter, striking a balance between global optimality and computational efficiency. Notably, the lower the threshold, the larger the subset of dimensions of the complicating variable, and vice versa. Its selection should align with the specific objective of the problem–opting for a higher threshold if fast feasible solutions are prioritized, and a lower threshold if global optimality is the primary goal.

| | **Algorithm 3: Learning to approximate complicating variables** |
|---|---|
| | **Inputs:** Mixed Integer Optimization model $MIP$, Trained $NN(\cdot)$ from Algorithm 2, target problem instance $\{\boldsymbol{A},\boldsymbol{B},\boldsymbol{b},\boldsymbol{c},\boldsymbol{d}\}^*$, probability threshold $k_{prob}$; |
| **1** | Approximate the complicating binary variable $y_m = NN(\{\boldsymbol{A},\boldsymbol{B},\boldsymbol{b},\boldsymbol{c},\boldsymbol{d}\}^*; k_{prob})$ and reduce dimension $\boldsymbol{M}$ to $\boldsymbol{M}^* \subseteq \boldsymbol{M}$ by obtaining the most likely predicted classes $m \in \mathbf{M}$; |
| **2** | Reduce the dimensionality of the original complicating variable to $y_{m^*} \in \{0,1\}, \ \forall\, m^* \in \boldsymbol{M}^* \subseteq \boldsymbol{M}$ |
| **3** | Solve $RMIP(\{\boldsymbol{A},\boldsymbol{B},\boldsymbol{b},\boldsymbol{c},\boldsymbol{d}\}^*)$ using a commercial solver; |
| | **Output**: Solution of $MIP(\{\boldsymbol{A},\boldsymbol{B},\boldsymbol{b},\boldsymbol{c},\boldsymbol{d}\}^*)$ |

# 2. Computational benchmark via an example case

The capabilities of the framework are assessed through its application to a mixed integer linear programming (MILP) model with a flow-capturing facility location allocation formulation that describes long-term investment planning and medium-term tactical planning in a personalized medicine supply chain for cell therapy manufacturing and distribution (Triantafyllou et al., 2022). The MILP model is presented in Appendix A.

### 3.1 Mixed integer optimization model

Chimeric Antigen Receptor (CAR) T-cell therapies are an advanced therapy medicinal product (ATMP). Owing to their promising clinical outcomes, they have paved the way for regulatory approval of numerous personalized cell therapy products (Young et al., 2022). Unlike typical cancer treatments, personalized therapies are characterized by patient-specific supply chains that consider bespoke parallel manufacturing lines and dedicated distribution nodes. The in-house mixed integer linear programming MILP model that describes the CAR T-cell therapy supply chain is used for the identification of optimal supply chain network structures (Triantafyllou et al., 2022). The MILP model is presented in Appendix A. As seen in Figure 2, the supply chain superstructure includes 4 nodes–leukapheresis centers, manufacturing sites, quality control (QC) sites, and hospitals. The CAR T-cell therapy supply chain begins at a clinical (leukapheresis) center, where T-cells are extracted from the patient's bloodstream, cryopreserved, and then sent to a manufacturing facility. During manufacturing, the leukapheresis sample undergoes a series of processing steps including genetic engineering, where the gene of interest that expresses the tumor-associated antigen receptor is inserted into the cells (Levine et al., 2017). Following this, the therapy undergoes in-house quality control (QC) to ensure product quality and safety before being shipped back to the hospital for administration to the patient. The model considers factors such as demand uncertainty, manufacturing capacity limitations, patient-specificity, and time and location constraints. The primary objective is to minimize the total supply chain cost, with the total turnaround time expressed as a non-monetary supply chain metric modeled as a constraint.

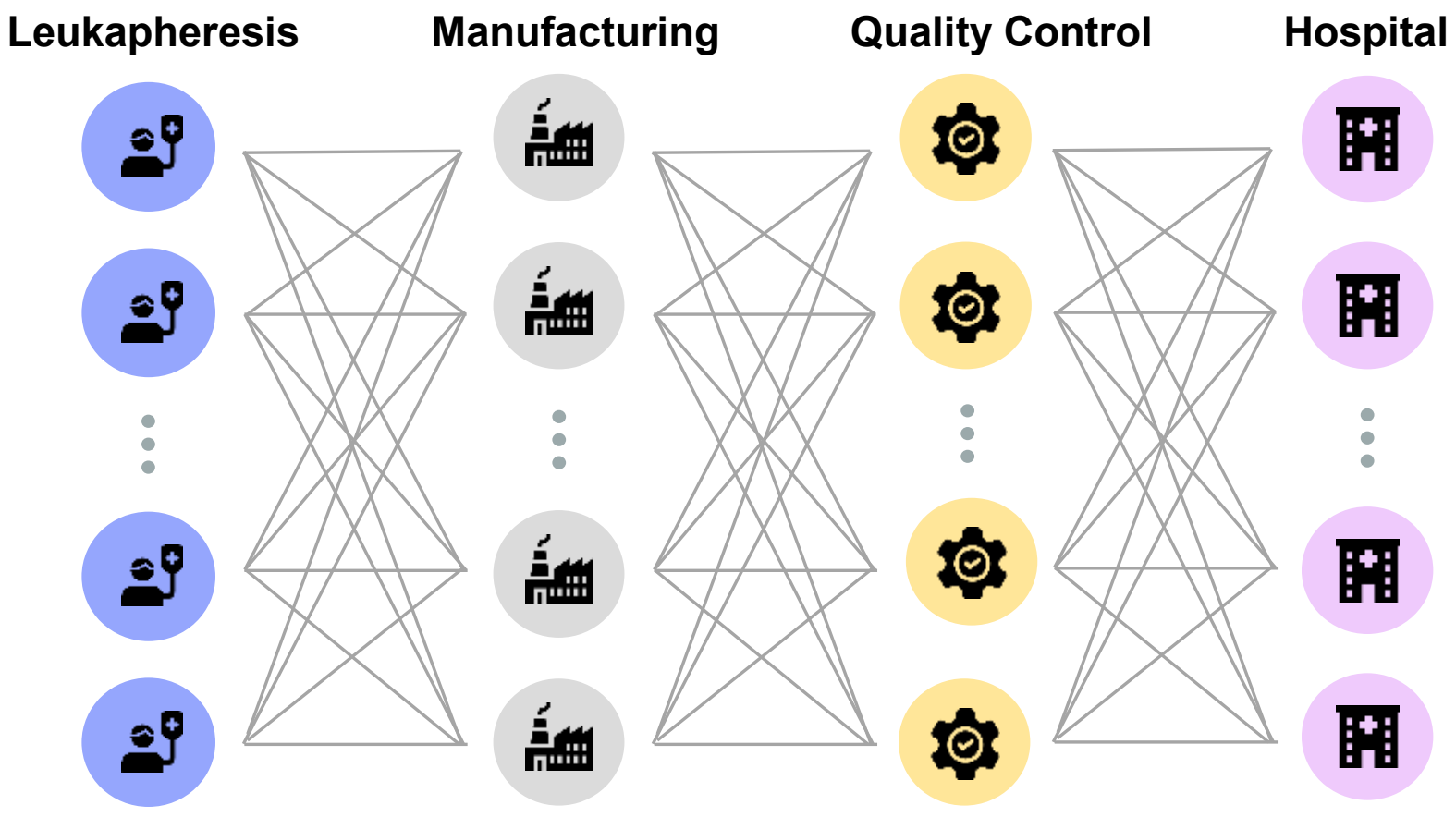

Figure 2. Representation of the CAR T-cell therapy supply chain superstructure. The main steps of a typical CAR T cell therapy lifecycle are (a) leukapheresis, (b) manufacturing, (c) Quality Control, and (d) therapy administration.

The supply chain manages the coordination and availability of different raw materials, alongside expert handling (e.g. cryopreservation) during the transportation of the samples and/or therapies. Ideally, all the processes involved should unfold within a short time frame (approximately 3 weeks). Concurrently, the patient-specific nature of the supply chain, which entails patient-specific manufacturing batches, places the patient schedule at the forefront (Papathanasiou et al., 2020). This implies that the allocation of available parallel lines in manufacturing facilities and the delivery process should be synchronized based on each patient's clinical condition and location. Presently, the CAR T-cell supply chain heavily relies on white-glove logistics and lacks a systematic approach to coordinated decision-making. In this regard, leveraging process systems engineering tools can guide decision-making by assisting manufacturers, suppliers, and clinicians in optimally coordinating these tasks for thousands of patients simultaneously.

Market forecasts anticipate the introduction of 60–90 new cell and gene therapy products, treating 60,000–90,000 patients annually by 2030 (Hillmoe and Shen, 2022; Quinn et al., 2019; Young et al., 2022). As demand grows, the supply chain optimization via MIP of this problem becomes challenging as its complexity scales exponentially with the increasing number of patients. Indicatively, for demands of up to 5000 therapies per year, the MILP model contains 68 million constraints and 16 million discrete variables. To address this, we utilize the proposed decomposition framework (Figure 1).

### 3.2 Complicating variable identification

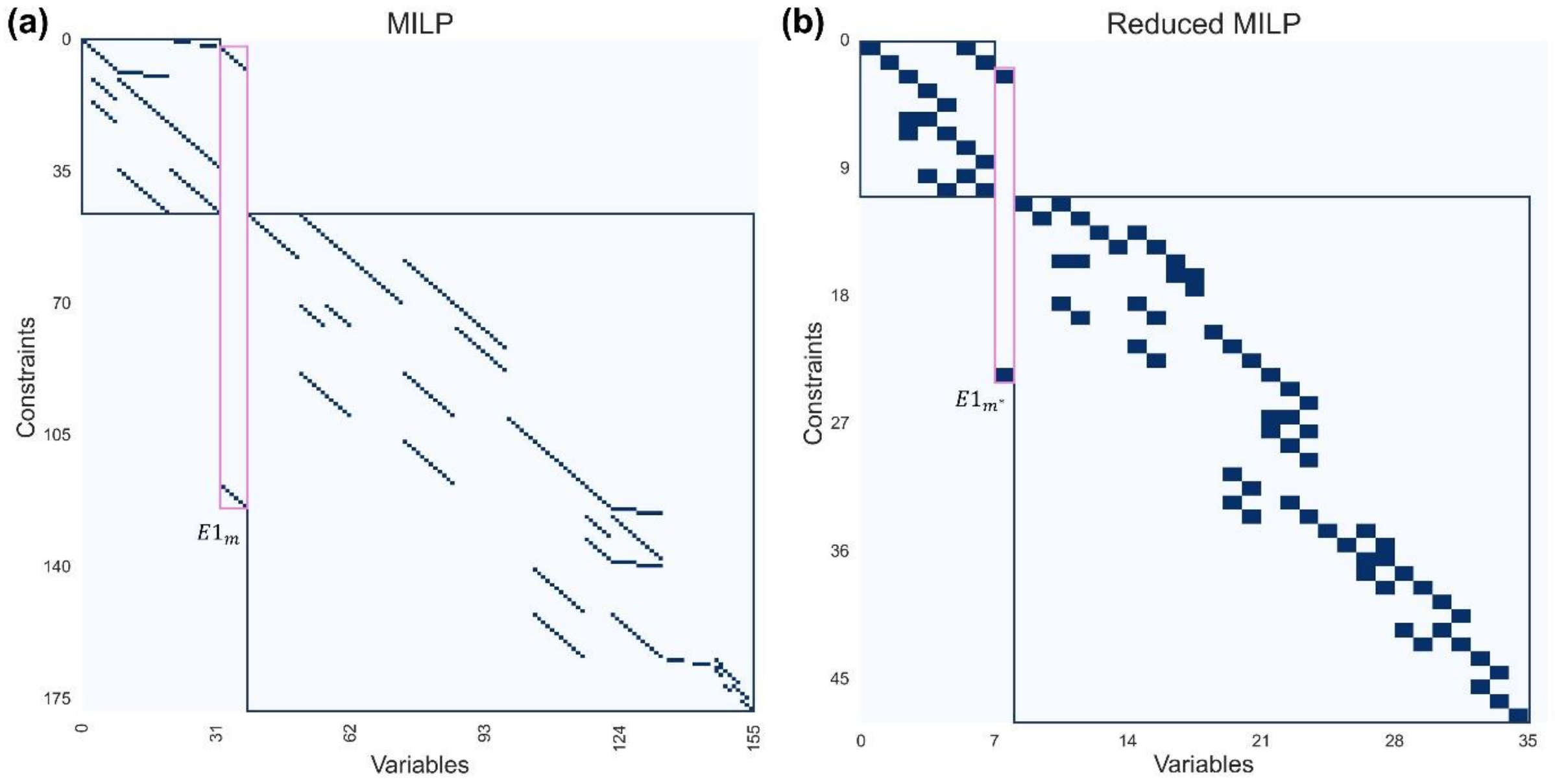

Figure 3. Constraint–variable matrix for a trivial instance in (a) the original MILP and (b) the reduced MILP. Complicating variables are highlighted with pink boxes, while decomposable subproblems are identified with blue boxes. The trivial instance includes a demand of 2 patients $p$, 1 leukapheresis center $c$, 6 candidate manufacturing sites $m$, 1 hospital $h$, 1 candidate transport mode $j$, and 1 time point $t$ (single-period model).

The first step of the proposed framework is to identify the complicating variables (Figure 3a) that will be approximated by the neural networks, which in this example case is the single-dimensional investment planning binary variable ($E1_m$) that determines the establishment of candidate manufacturing facilities $m$. If this variable is fixed, the problem can be easily solved by blocks as it will be separated in a strategic planning (facility location) problem and a medium-term scheduling (patient allocation) problem. By approximating the investment planning binary variable, the manufacturing facilities chosen for each problem instance can be determined by the deep learning model, allowing the MILP to be solved with a reduced set of facilities. The solver can then select some or all of the manufacturing facilities predicted by the deep learning model by solving the reduced MILP, aiming to achieve the optimal solution. The resulting reduced MILP is a more manageable model, as demonstrated by the substantial reduction in constraints and variables illustrated in Figure 3b for a trivial problem instance.

### 3.3 Dataset generation

The data for the nominal problem instance are obtained from the literature (Triantafyllou et al., 2022). We consider 4 leukapheresis sites, 4 hospitals in the UK, and 6 manufacturing sites located in the UK, US, and Europe. The manufacturing facilities $m$ have a capacity of 4 ($m_1$ and $m_4$), 10 ($m_3$ and $m_6$), or 31 ($m_2$ and $m_5$) parallel lines, and a forward-looking scenario of 8 days of manufacturing is considered.

Current market size and the demand for these products are highly uncertain (Papathanasiou et al., 2020), therefore, the uncertain parameter that varies throughout the instances is the three-dimensional demand profile $INC_{p,c,t}$ indicating which patient $p$ walks to which leukapheresis center $c$ at which day $t$. We assume that the demand profile is repeated every trimester and we generate random instances for $INC_{p,c,t}$ (Algorithm 1). For this problem case, we consider demands ranging from 8 to 2400 patients per year (or 2 to 600 patients per trimester) with different probability distributions. The probability distributions used are uniform, left triangular, and right triangular (Figure 2), based on discussions with experts. For each total demand scenario (number of patients per trimester) and each distribution, 5 randomized demand profiles are generated. In total, $3(distributions) \times 5(\frac{profiles}{distribution}) \times (600 - 2) = 8{,}970$ instances are generated. All generated demand profiles are bound to satisfy the feasibility constraint that a maximum of 8 patients per day can visit a leukapheresis center.

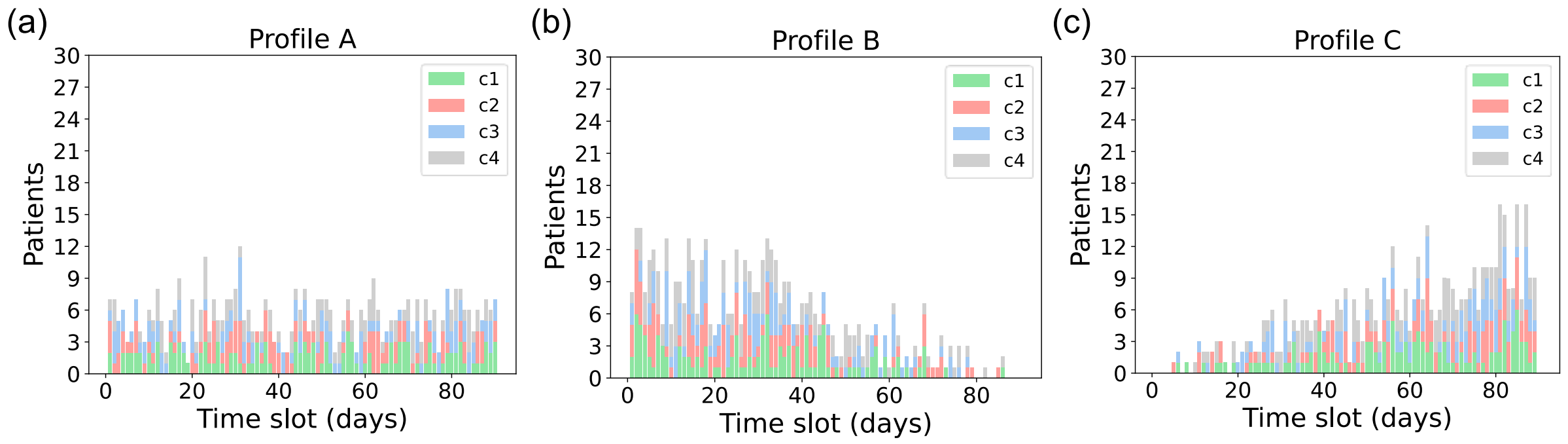

Figure 4. Randomized demand profiles ($INC_{p,c,t}$) for the 2000 patients' demand scenario for 3 different probability distributions: (a) uniform distribution (b) left triangular distribution (c) right triangular distribution (Triantafyllou et al., 2023).

To generate the dataset for training the multi-label deep learning classifier, we employ 7 labels; namely: the 6 manufacturing facilities and the infeasibility label that indicates when a problem instance is infeasible. The labels for each instance are derived by solving the MILP for all 8,970 problem instances with the IBM ILOG CPLEX 20.1 optimizer to global optimality (`mipgap= 0.0001`) in Pyomo 6.7 (Hart et al., 2011). We perform this using high-performance computing (HPC), parallel computing through multiple CPUs, and array jobs to expedite the synthetic data generation process (24 parallel 2x Intel Icelake Xeon Platinum 8358 cores and a maximum RAM of 64GB ("Imperial College Research Computing Service," 2023). Each instance took 18 to 4,186 CPU seconds to run.

Figure 5 illustrates the label distribution within the dataset, revealing a highly imbalanced dataset, where $m_2$ is preferred over other facilities. This is inherent to the nature of the problem, where demand zones (i.e., hospitals) are exclusively located in the UK around London. As a result, manufacturing facilities closer to these demand zones always receive higher priority. Notably, among the two large facilities $m_2$ and $m_5$, $m_2$ is always prioritized despite being located in Europe. This prioritization is because Belgium (EU) is closer to the demand zones compared to Glasgow (UK). This imbalance is a significant challenge when training the classifiers, as it reduces the ability to effectively learn from and predict underrepresented classes.

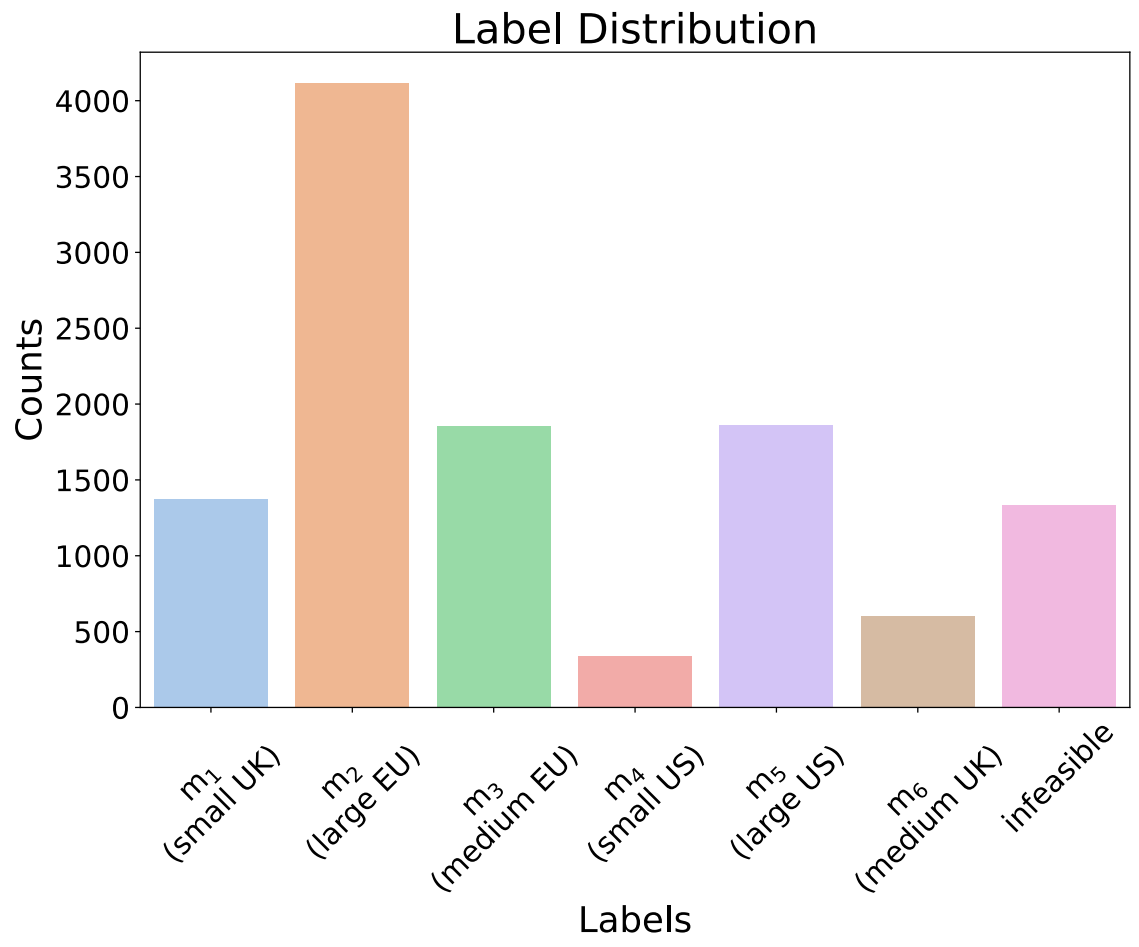

Figure 5. Label distribution for the benchmark MILP problem.

## 3.4 Deep learning models

The dataset is split into three subsets: 80% for training, 15% for testing, and 5% for validation. We consider and compare their performance in two different architectures: (a) feed-forward neural networks (ANN) and (b) convolutional neural networks (CNN). The input features considered for the networks are the total daily demands for a quarter of a year (90 days) ($INC_{p,c,t}$), assuming a recurrent demand profile per trimester. The output layer includes 7 neurons, where the labels are the 6 manufacturing facilities and the possibility of infeasible solutions due to limited facility capacity (Figure 6). Thus, each network consistently processes 90 input features, regardless of the significant variations in the number of patients and therefore the variation in the number of constraints and variables.

Algorithm 2 is employed to tune the hyperparameters and train the model for both architectures. The Adam optimizer Binary Cross Entropy with Logits loss function is used for both architectures. We set $maxiter = 300$ and perform Bayesian optimization with the sample level accuracy as the objective function that is maximized. We implement the neural network models in PyTorch(Paszke et al., 2019), and use GPyTorch (Gardner et al., 2018) to build the gaussian process surrogates. BoTorch(Balandat et al., 2019) is used to perform Bayesian optimization. The neural network training is performed on GPU, significantly reducing the training time. Bayesian optimization is performed on the CPU. The GPU NVIDIA Quadro RTX 6000 is used ("Imperial College Research Computing Service," 2023). Notably, the training of 300 neural networks with Bayesian optimization requires 4.91 h and 56.48 h for ANNs and CNNs, respectively. As the architecture of the neural networks remains unaffected by the size of the problem (number of constraints and variables), training times are expected to remain in the same orders of magnitude even if the size of the MIP instances substantially increases.

The performance of Algorithm 2 is presented in Figure 7. It is evident that the CNN converges to the optimal hyperparameters more rapidly than the ANN. This can be attributed primarily to the CNN's hyperparameters being constrained to non-architectural configuration ones.

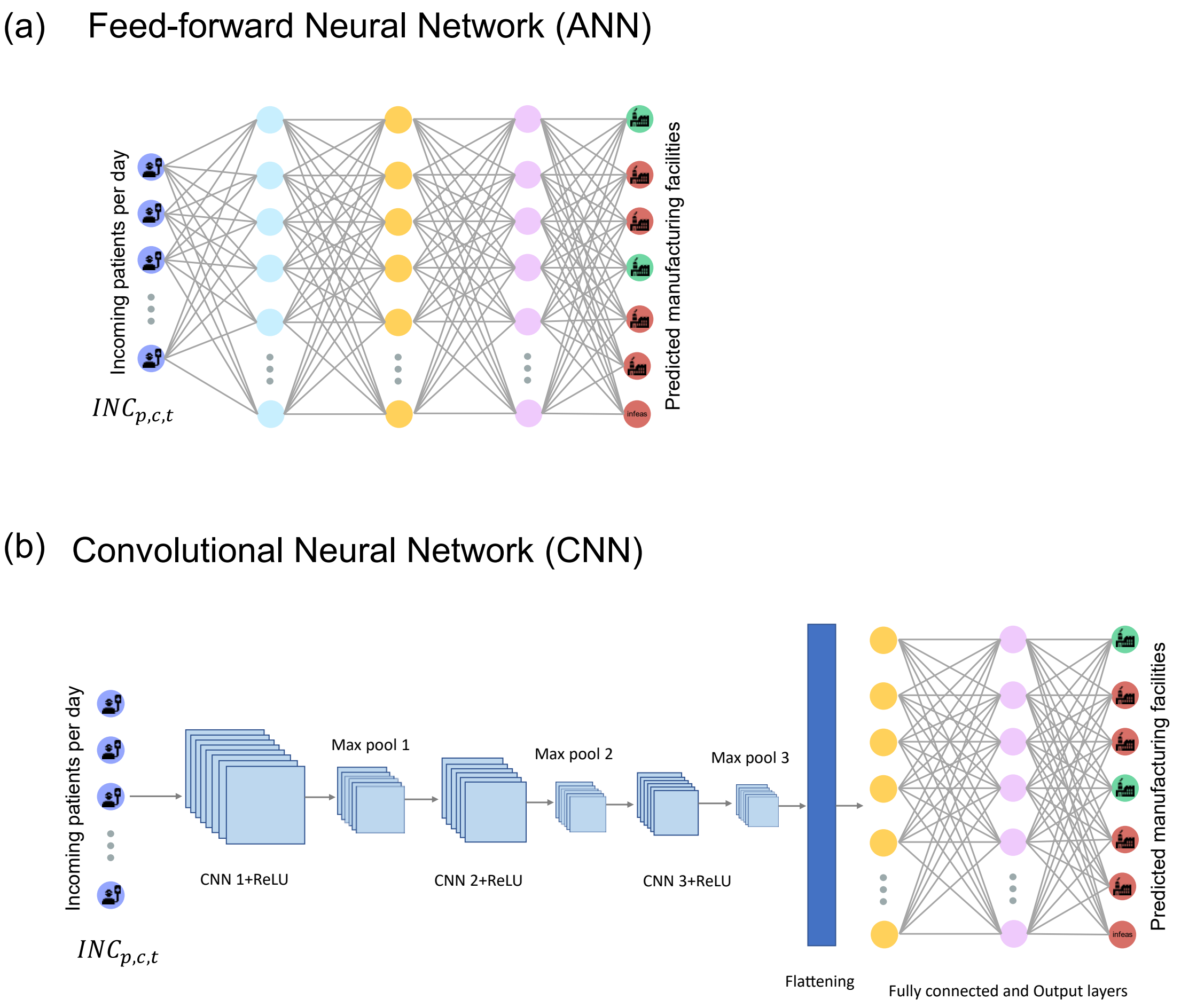

Figure 6. Neural network configuration for (a) feedforward neural networks and (b) convolutional neural networks.

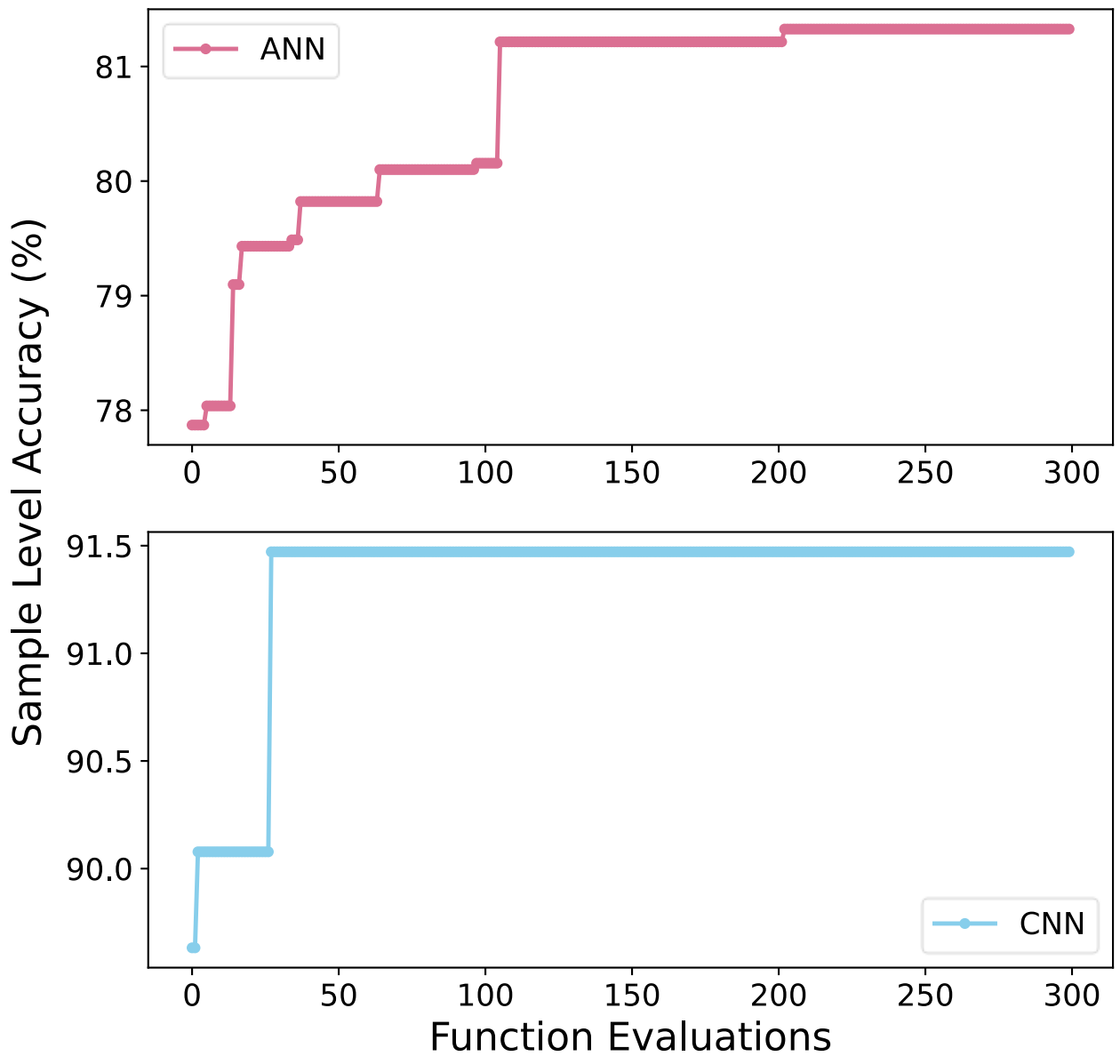

Figure 7. Bayesian optimization progress for hyperparameter tuning for feedforward neural networks and convolutional neural networks.

### (a) Feed-forward neural networks

We aim to maximize the multi-label classifiers' performance towards global optimality and we tune the: number of hidden layers ([1, 3]), number of neurons in the hidden layers ([50, 256]), learning rate on a logarithmic scale ([$10^{-5}$, $10^{-1}$]), number of training epochs ([500, 15000]). The final configuration of the hyperparameters is shown in Table 1.

Table 1. Hyperparameter configuration of the feed-forward neural network multi-label classifier.

| **Hyperparameters** | **Feed-forward neural network** |
|---|---|
| Epochs | 8715 |
| Optimizer | Adam |
| Activation function | ReLU |
| Learning rate | $6.19 \cdot 10^{-5}$ |
| Hidden layers | 3 |
| Neurons in the hidden layers | 200 |

### (b) Convolutional neural networks

For the convolutional neural networks, the learning rate on a logarithmic scale ([$10^{-5}$, $10^{-1}$]), number of training epochs ([500, 15000]), dropout rate in the 1st dropout layer ([0, 0.3]), and dropout rate in the 2nd dropout layer ([0, 0.3]) are tuned. To cater for the interdependency of the hyperparameters in CNN architectures which may challenge the BO, here, the layers and kernel sizes were tuned manually. The final configuration is given in Table 2.

Table 2. Hyperparameter configuration of the convolutional neural network multi-label classifier.

| **Hyperparameters** | | **Convolutional neural network** |
|---|---|---|
| | Epochs | 7532 |
| | Optimizer | Adam |
| | Activation function | ReLU |
| | Learning rate | 4.37·10-4 |
| Convolution layer 1 | Input channels | 32 |
| | Output channels | 64 |
| | Kernel size | 10 |
| | Stride | 1 |
| | Padding | 1 |
| Max pooling layer 1 | Kernel size | 5 |
| | Stride | 5 |
| | Padding | 0 |
| Convolution layer 2 | Input channels | 64 |
| | Output channels | 128 |
| | Kernel size | 5 |
| | Stride | 1 |
| | Padding | 1 |
| Max pooling layer 2 | Kernel size | 3 |
| | Stride | 3 |
| | Padding | 0 |
| Convolution layer 3 | Input channels | 128 |
| | Output channels | 256 |
| | Kernel size | 3 |
| | Stride | 1 |
| | Padding | 1 |
| Max pooling layer 3 | Kernel size | 2 |
| | Stride | 3 |
| | Padding | 0 |
| Fully connected layer 1 | Input neurons | 256 |
| | Output neurons | 256 |
| | Dropout rate | 0.0 |
| Fully connected layer 2 | Neurons | 256 |
| | Output neurons | 128 |
| | Dropout rate | 0.3 |

### 3.5 Multi-label classifier evaluation

#### 3.5.1 Performance metrics

The overall performance of the two deep learning classifiers is evaluated through the metrics summarized in Table 3. Further details on how those are calculated are provided in Appendix B (Tables B1 and B2). The normalized multi-label confusion matrix is presented in Figure 8 for both neural network architectures.

Table 3. Evaluation metrics for the feed-forward neural network.

| | **Feed-forward neural network** | **Convolutional neural network** |
|---|---|---|
| Sample-level accuracy (%) | 81.32 | 91.47 |
| Jaccard Index (%) | 94.5 | 97.76 |
| Hamming Loss (%) | 4.96 | 2.23 |

CNNs demonstrate an improved classifier performance, with sample-level accuracy, Jaccard Index, and Hamming Loss outperforming its feedforward neural network counterpart. Notably, the convolutional architecture predicts the global optimum solution in 92.47% of instances within the test set, demonstrating its efficacy in capturing intricate patterns and temporal relationships in the demand profiles.

Classes $m_2$, $m_4$, $m_5$, and the infeasible class, all exhibit high precision, recall, and F1-scores, indicating robust performance for both neural network architectures (Tables B1 and B2). Class $m_6$ shows lower recall than precision, indicating that true positive instances (TP) are not always captured. Classes $m_1$ and $m_3$ demonstrate balanced precision and recall values in the ANN, resulting in moderate F1-scores, while class $m_3$ in the CNN shows lower recall than precision.

#### 3.5.2 Confusion matrices

By examining the confusion matrix in Figure 8a, it can be seen that, in 26% of the instances, the model confuses facility $m_1$ (4 parallel lines) with facility $m_3$ and $m_6$ (10 parallel lines) and with the no-label prediction. Specifically, the model accurately predicts that $m_1$ belongs to the optimal supply chain configuration in 74% of the cases. However, it misclassifies facility $m_1$ with $m_3$ in 8% of the cases, or it misclassifies facility $m_1$ with $m_6$ in 3% of the cases. In 13% of the cases, it fails to predict $m_1$ as an extra required facility. These misclassifications can be attributed to the model tending to overpredict by choosing a slightly bigger facility to meet the demand and in other cases, it underpredicts by omitting an extra facility. In the former scenario, the deep learning model guides the MILP to a local optimum solution, whereas the latter case leads to infeasible solutions. In contrast, the CNN demonstrates improved performance, accurately predicting $m_1$ in 89% of the cases and misclassifying it less with $m_3$ and $m_6$ (4% of occurrences), while it fails to predict $m_1$ in 7% of the cases. A similar behavior is observed with facility $m_6$ (10 parallel lines). It should be noted that the classifier is performing well when predicting infeasibilities, reaching a precision of 94% for the ANN and 97% for the CNN. On rare occasions, the classifier considers infeasible solutions as feasible. In those cases, the reduced MILP will determine the infeasibility of the solutions.

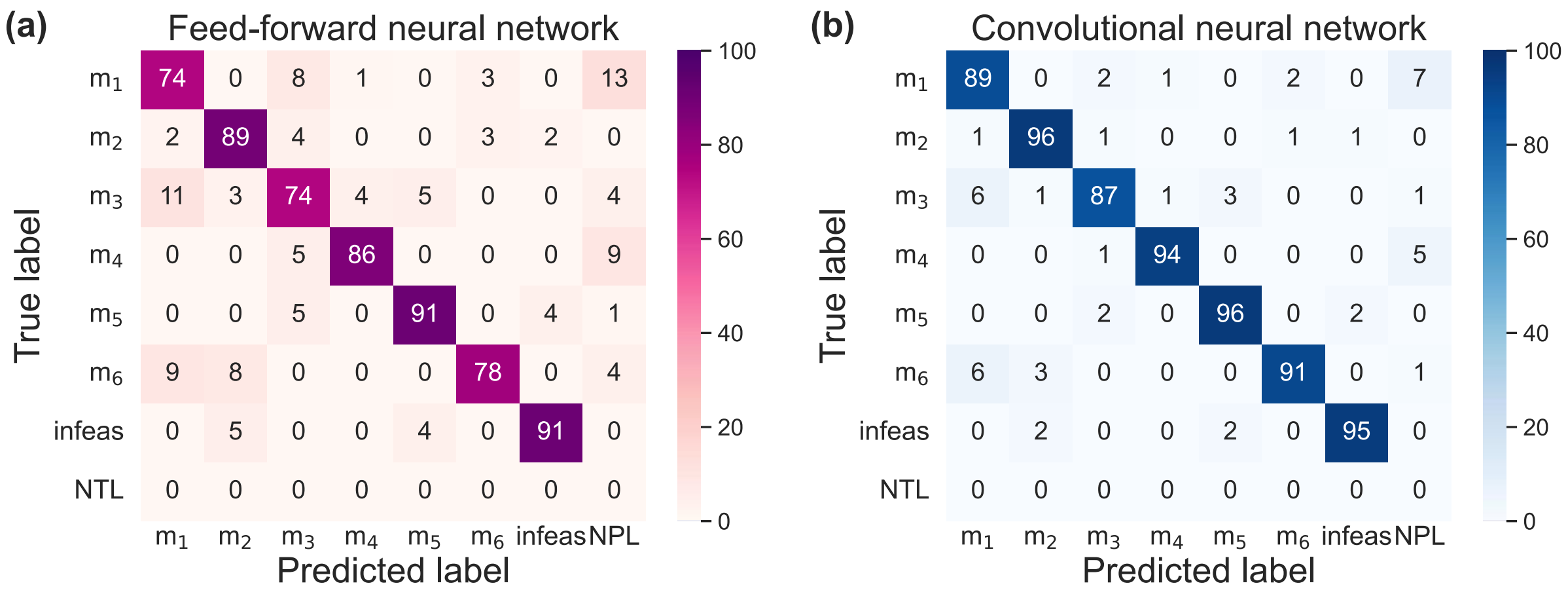

Figure 8. Normalized multi-label confusion matrix with one extra row for No True Labels (NTL) and one extra column for No Predicted Labels (NPL) for (a) the feedforward neural networks and (b) the convolutional neural networks.

### 3.5.3 Deep learning model complexity – accuracy trade-off

Training CNNs is computationally more expensive than training ANNs, especially as the model complexity increases. Each evaluation during hyperparameter tuning involves training the model, and the computational cost can quickly escalate with CNNs, leading to longer optimization times. The hyperparameter tuning process for the CNN took 55 hours and 37 minutes for 300 function evaluations, significantly surpassing the time required for 300 ANN trainings, which took 4 hours and 57 minutes. However, as shown in Figure 7, the CNN converged to the optimum hyperparameters within 28 iterations, totaling a computational time of 5 hours and 21 minutes. In contrast, the ANN achieved convergence at iteration 203, taking 3 hours and 19 minutes. The efficiency of the CNN's convergence is attributed to the simpler nature of its hyperparameters. Unlike the ANN, the CNN does not involve the optimization of network architecture specifics, such as convolutional, max pooling, fully connected, and dropout layers.

### 3.6 Reduced mixed integer optimization model

In section 3.4 we demonstrated the development of the neural networks that are responsible for long-term strategic planning, predicting the correct manufacturing facilities for the target MILP instance. Here, the solution of the NNs is used to reduce the MILP model that is solved for a subset of manufacturing facilities as chosen by the classifier. The reduced MILP formulation is treated as a subproblem of the original MILP model and addresses medium-term tactical planning in the supply chain. The objective is the minimization of the therapy cost and return time by identifying good candidate solutions for the: allocation of patient samples in the manufacturing facilities and hospitals, transport modes for the node-to-node connections, and the utilization of the available parallel lines in the manufacturing sites. It should be noted, however, that similar to many heuristic or decomposition approaches, the algorithm does not provide a guarantee of global optimality.

Figures 9 and 10 illustrate the performance of the deep learning decomposition framework for 12 instances from the validation set with different probability distributions for the $INC_{p,c,t}$ parameter (100, 200, 500, 1000, and 2000 patients annually) for the feed-forward and convolutional neural networks, respectively. The global optimal solutions of the MILP and the predictions of the ANN and CNN classifiers are provided in Appendix C (Tables C1 and C2). For uniform demand distributions (Figures 9a, 9d, 10a, and 10d), the decomposition framework attains the global optimum solution for all instances for both neural network architectures. In the ANN, there is a notable reduction in constraints and binary variables reaching up to 66.6% and 64.7% respectively. The CNN outperforms the ANN in the 1000 patients per year instance (Figures 10a and 10d). In this scenario, both the 50% and 1% probability thresholds lead to identical solutions for the CNN, predicting that only facility $m_2$ should be built. This is because the CNN is more confident, accurately predicting the correct manufacturing facility and no extra facilities reducing the constraints by 86.9% and the binary variables by 83.3%.

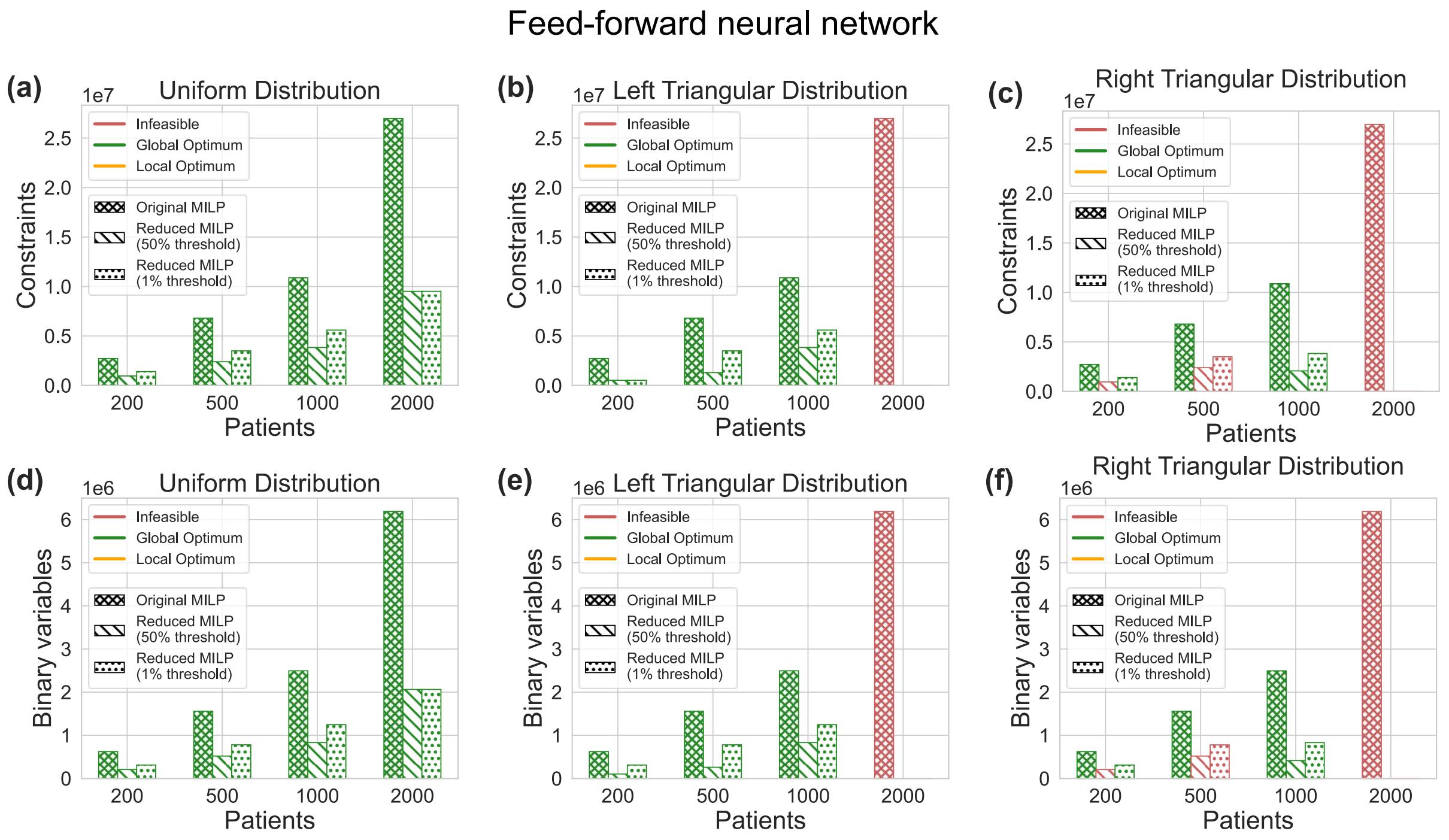

Figure 9. Performance of the feed-forward neural network decomposition framework in terms of constraints and binary variables in comparison to the original MILP model for different probability distributions of 12 parameter instances from the validation set. Two distinct scenarios are examined: one where the neural network employs a 50% probability threshold for predictions, and another scenario where a more conservative 1% threshold is utilized.

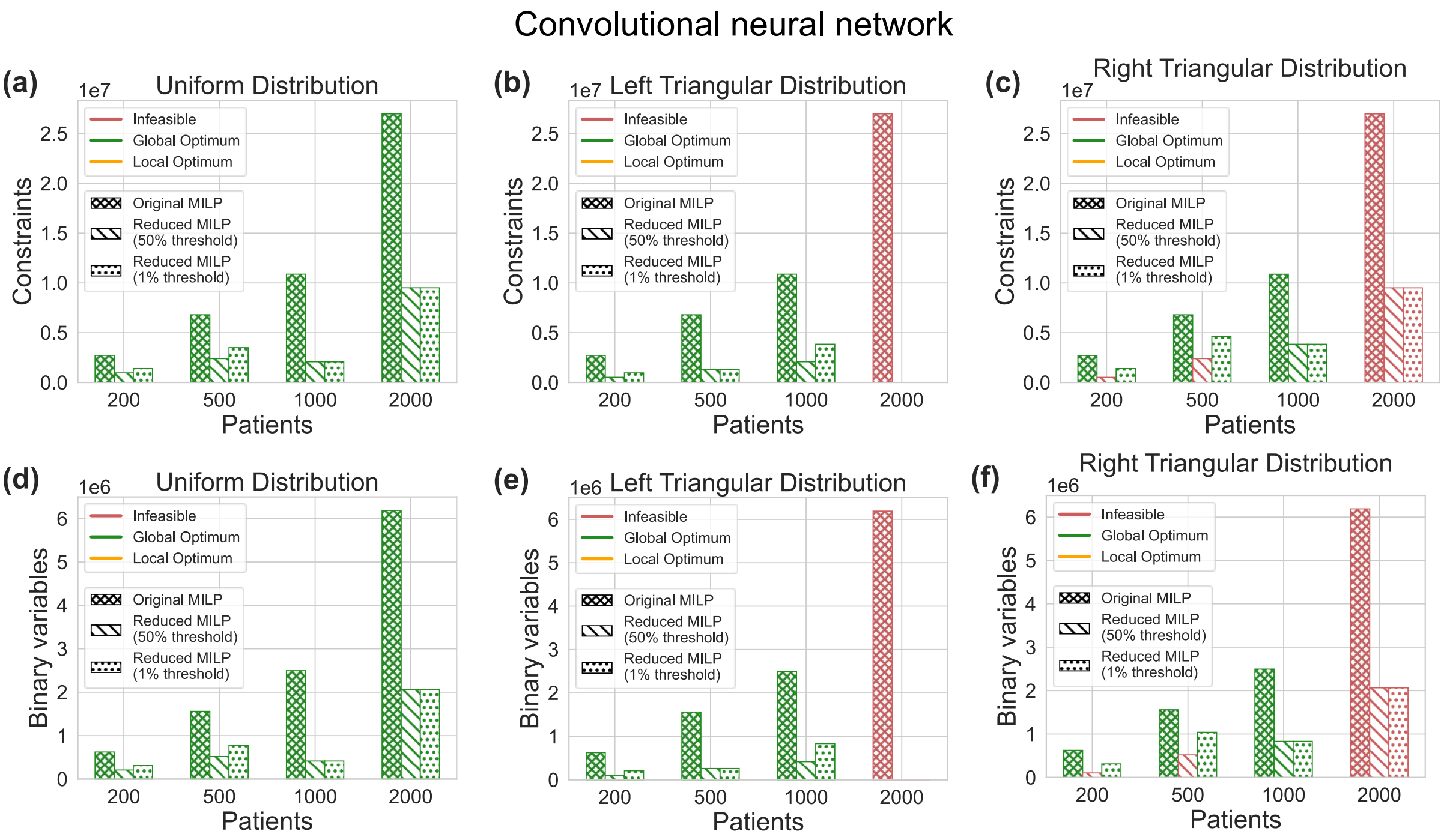

Figure 10. Performance of the convolutional neural network decomposition framework in terms of constraints and binary variables in comparison to the original MILP model for different probability distributions of 12 parameter instances from the validation set. Two distinct scenarios are examined: one where the neural network employs a 50% probability threshold for predictions, and another scenario where a more conservative 1% threshold is utilized.

In the case of the left triangular distribution (Figures 9b, 9e, 10b, and 10e), the decomposition framework identifies the global optimum solutions across all instances. Particularly noteworthy is the case of 2,000 patients per year, where the classifier deems the problem infeasible, resulting in no solution attempts and consequently reducing the overall computational time. Overall, the CNN achieves faster solutions by predicting fewer dimensions (manufacturing facilities). This is highlighted in the 1% threshold predictions, where the manufacturing facilities that have a probability of more than 1% being optimal are less compared to the ANN.

In the right triangular distribution (Figures 9c, 9f, 10c, and 10f), the decomposition framework generates infeasible solutions for the 200 and 500 patient demand scenarios per year for the 50% threshold, despite the actual feasibility of the problems. This discrepancy can be attributed to underpredictions, where smaller facilities than required are assigned. The classifier with the 1% threshold performs better, guiding the solver to the global optimum solution for all instances, including both 200 and 500 patient demand scenarios in the CNN. In the ANN, even with the 1% probability threshold, the 200-patient demand scenario is wrongfully deemed infeasible. For the instance of 2,000 patients per year, the ANN outperforms the CNN, by predicting the infeasible instance as infeasible. The CNN classifier predicts feasibility by predicting facilities $m_2$ and $m_5$ leading to a 64.6% reduction in the number of constraints and a 66.6% reduction in the number of binary variables. However, upon solving the reduced model, the solver deems the model infeasible (Figures 10c and 10f).

There are no local optimum solutions in the 12 instances depicted in Figures 9 and 10. This is owing to the fact that the neural networks make predictions, but then instead of directly fixing the binary variable based on these predictions, the dimensions are reduced. This reduction ensures that the MILP dictates the optimal facilities from the subset provided to it. By allowing the MILP model to make the final decisions within the reduced subset, the potential for local solutions is effectively mitigated.

# 3. Conclusions and outlook

Industrially relevant Mixed-Integer Problems (MIPs) are often complicated and computationally expensive to solve. In this context, Machine Learning (ML)-enhanced methodologies offer a route to decrease the problem complexity and enable the solution of larger formulations by exploiting problem-specific and distribution-specific structures. While many commercial solvers are designed to handle a broad range of problems effectively, they may limit user accessibility to tailor them to specific problems. In contrast, learning-based approaches can adapt solvers to particular applications that are of interest to the user. By training an ML model on instances drawn from specific problem distributions, the ML-enhanced optimizer performance improves and can outperform generic solvers that do not account for such characteristics. This approach is particularly useful in situations where multiple related instances need to be solved, and while the constraints and variables of each instance change significantly, they do so in a predictable manner.

In this work, we present a 'learning to reduce model dimensionality' approach to address the computational complexity of large-scale MIP models, by learning to construct problem-specific heuristics from a dataset of wide MIP instances that belong to different distributions. For the development of the presented approach, we train and compare two multi-label deep learning classifiers. The first one is using Artificial Neural Networks (ANNs), while the second one is based on Convolutional Neural Networks (CNNs). In both cases, the objective of the multi-label classifier trained on precomputed MIP instances is to identify the active dimensions of the MIP. This is achieved by approximating (a) complicating variable(s). We introduce a methodology for generating synthetic MIP instances that addresses parameter uncertainty, assisting the deep learning classifiers in generalizing better to instances from multiple distributions. Critically, we additionally introduce an "infeasible" label that enables the classifier to predict when an instance is infeasible. This is highly important as it can reduce the amount of computational time one will need to build and solve infeasible large-scale MIP instances. Once the active dimensions are identified by the classifier, the reduced MIP is solved for the subset of dimensions with a standard off-the-shelf solver. Our methodology is driven by a user-specified parameter (probability threshold) that allows to directly trade-off the conservativeness of the predictions. We test the performance of the algorithms on a large-scale, complex case of personalized medicine supply chain, where we approximate investment planning binary variables. In this context, both the ANN and CNN classifiers perform admirably. Notably, the CNN classifier leads the reduced model to the global optimum solution for 91.47% of the instances in the test set, while the ANN does it for 81.32% of the instances. With respect to their performance in terms of accurately predicting the infeasibility of the solution, both classifiers score highly with a $precision \geq 94\%$. Both the ANN and the CNN classifier succeed in significantly decreasing the number of constraints by a maximum of 81% and 83% respectively.

To further enhance the presented approach by identifying better problem partitions, one could explore the integration of automated complicating variable identification and graph visualization techniques (Allman et al., 2019; Jalving et al., 2020). Additionally, learning-based decomposition approaches as proposed by Mitrai and Daoutidis, (2024b) could be explored to assess the feasibility of further decomposition, along with initialization strategies (Mitrai and Daoutidis, 2024a). Future work could involve the application of learning strategies within the Branch-and-Cut, such as developing problem-specific heuristics to optimize variable branching (Nair et al., 2020), node selection(Labassi et al., 2022), and cutting plane selection (Li et al., 2023) tailored to specific problems.

Concluding, ML-enhanced methods are reliant on the solution of the MIP instances for their one-off training, which may render them computationally expensive. Notwithstanding the foregoing, this trade-off may be tolerated particularly for large-scale, real-world problems, as the generated classifier will thereafter enable the expedited solution of MIPs, compensating for the initial computational expense. It is important to highlight that our framework efficiently handles scalability, as it is agnostic of the exact constraints and variables of the problem ($\boldsymbol{Ay} + \boldsymbol{Bx} \leq \boldsymbol{b}$) and only handles parameters ($\{\boldsymbol{A}, \boldsymbol{B}, \boldsymbol{b}, \boldsymbol{c}, \boldsymbol{d}\}$). The architecture of the neural networks remains unaffected by the size of the problem (number of constraints and variables), as the number of inputs and outputs are consistent across different MIP instances. This choice of input features ensures that both scalability and training times are not adversely impacted as the problem size increases. However, other case studies may involve more complex dynamics with multi-dimensional uncertain parameters that cannot be simplified. In such instances, architectures like graph neural networks (Scarselli et al., 2009) could provide a more detailed representation. These networks are capable of encoding all constraints and variables of the problem (Mitrai and Daoutidis, 2024b), potentially increasing predictive accuracy at the expense of increased computational time.

**Acknowledgments**

Funding from the UK Engineering & Physical Sciences Research Council (EPSRC) for the Future Targeted Healthcare Manufacturing Hub hosted at University College London with UK university partners is gratefully acknowledged (Grant Reference: EP/P006485/1). Financial and in-kind support from the consortium of industrial users and sector organizations is also acknowledged. MMP and NT would like to thank Matthew Lakelin (Trakcel Ltd) for the expert input provided in the development of the previously published cell therapy supply chain model. NT is thankful for the Marit Mohn Scholarship awarded by the Department of Chemical Engineering, Imperial College London.

**Data availability**

The codes and dataset developed in this study are available open-source.

## Appendix A

The MILP model for the CAR T-cell therapy supply chain is presented here (Triantafyllou et al., 2022). It is based on the following assumptions:

- The leukapheresis sites are pre-defined, aligned with the UK healthcare system
- Each therapy needs to return to the hospital that is co-located with the leukapheresis site from which the original sample was taken
- The capacity of the candidate manufacturing facilities is fixed
- The demand profile is repeated every trimester of the year

The objective of the model is to minimize the total cumulative cost of all manufactured therapies (Eq. 1).

$$min\, TOTCOST = \sum_p CTM_p + \sum_p TTC_p + \sum_p CQC \qquad (1)$$

The total cost (Eq. 1) is calculated as the sum of the manufacturing ($CTM_p$), transport ($TTC_p$) and quality control costs ($CQC$) for all patients $p$. In the scenarios presented here, quality control is considered an in-house procedure with a fixed equal cost for each patient.

*Manufacturing costs.* The first term of the objective function ($\sum_p CTM_p$) represents the manufacturing cost of all therapies $p$. Equation (2) calculates the manufacturing cost per therapy $p$.

$$CTM_p = \frac{\sum_m (E1_m \cdot (CIM_m + CFVM_m))}{NP} + CVM_p \,, \ \forall\, p \quad (2)$$

The first term of Equation (2) ($\frac{\sum_m E1_m \cdot CIM_m}{NP}$) corresponds to the capital investment for the construction of new facility $m$ as attributed to each therapy $p$. The second term ($\frac{\sum_m E1_m \cdot CFVM_m}{NP}$) reflects the fraction of variable manufacturing costs, also referred to as "fixed-variable costs", that cannot be easily adjusted to variable demand once the target productivity has been decided. Those may include personnel, facilities and equipment maintenance costs. Such costs are ~80% of the total variable costs and are only dependent on the number of established parallel lines, regardless of the utilization of the facility, therefore they can lead to a significant increase of the average cost per therapy if the average utilization of a manufacturing facility is low (Spink and Steinsapir, 2018) . The third term ($CVM_p$) describes the variable costs associated with the materials required for the manufacturing. Given the autologous nature of the CAR T therapies, this cost is only marginally affected by economies of scale and we assume that it is independent of the size of the facility.

Eq. (3) calculates the percentage of utilization of facility $m$ at time $t$.

$$RATIO_{m,t} = \frac{\sum_p DURM_{p,m,t}}{FCAP_m} \,, \ \forall\, m, t \quad (3)$$

$DURM_{p,m,t}$ is a continuous variable that takes a value of 1 only for the time period $t$ in which a therapy $p$ is being manufactured in facility $m$.

*Transport costs.* The second term in Eq. (1) ($\sum_p TTC_p$) corresponds to the total cost of transport for all therapies $p$, while the transport cost per therapy $p$ ($TTC_p$) is given by Eq. (4). The two terms correspond to transport from: (1) clinical site to manufacturing and (2) manufacturing site to hospital, respectively.

$$TTC_p = \sum_{c,m,j,t} Y1_{p,c,m,j,t} \cdot TT1_j \cdot U1_{c,m,j} + \sum_{m,h,j,t} Y2_{p,m,h,j,t} \cdot TT2_j \cdot U2_{m,h,j} \,, \ \forall\, p \quad (4)$$

*Material balances*. The model formulation includes a series of material balances around the facilities that are described by Eq. (5)-(12). Patient samples $p$ collected at the leukapheresis site $c$ at time $t$ ($INC_{p,c,t}$) will be ready to be shipped to the manufacturing facility ($OUTC_{p,c,t}$ ) after the duration of the leukapheresis procedure ($TLS$):

$$INC_{p,c,t} = OUTC_{p,c,t+TLS}, \forall p,c,t \tag{5}$$

Patient samples $p$ collected at the leukapheresis site $c$ being shipped to manufacturing facility $m$ via transport mode $j$ at time $t$ ($LSR_{p,c,m,j,t}$) will arrive at the manufacturing facility ($LSA_{p,c,m,j,t}$) after the duration of the transport activity ($TT1_j$) (Eq. 6):

$$LSR_{p,c,m,j,t} = LSA_{p,c,m,j,t+TT1_j}, \forall p,c,m,j,t \tag{6}$$

For each patient sample $p$, the outgoing samples $OUTC_{p,c,t}$ from a leukapheresis site $c$ at time $t$ are equal to patient sample $p$ sent to all manufacturing facilities $m$ under any transport mode $j$ at time $t$ ($\sum_{m,j} LSR_{p,c,m,j,t}$) (Eq. 7).

$$OUTC_{p,c,t} = \sum_{m,j} LSR_{p,c,m,j,t}, \forall p,c,t \tag{7}$$

Patient sample $p$ entering manufacturing facility $m$ at time $t$ is equal to patient samples $p$ shipped from all leukapheresis sites $c$ to manufacturing facility $m$ under any transport mode $j$ ($\sum_{c,j} LSA_{p,c,m,j,t}$) (Eq. 8).

$$INM_{p,m,t} = \sum_{c,j} LSA_{p,c,m,j,t}, \forall p,m,t \tag{8}$$

The outlet of patient samples $p$ of a manufacturing facility $m$ at time $t$ ($OUTM_{p,m,t}$) will be ready for shipment after the manufacturing process ($TMFE$) and quality control ($TQC$) have been completed (Eq. 9).

$$INM_{p,m,t} = OUTM_{p,m,t+TMFE+TQC}, \forall p,m,t \tag{9}$$

The sample $p$ leaving the manufacturing facility $m$, $OUTM_{p,m,t}$, is equal to the sample ready to be shipped to any hospital $h$ ($MSO_{p,m,h,j,t}$) under transport mode $j$:

$$OUTM_{p,m,t} = \sum_{h,j} MSO_{p,m,h,j,t}, \forall p,m,t \tag{10}$$

$MSO_{p,m,h,j,t}$ sample $p$ that has left manufacturing facility $m$ arrives at hospital $h$ via transport mode $j$ at time $t$ ($FTD_{p,m,h,j,t}$) after $TT2_j$ that is the duration of transport via transport mode $j$:

$$FTD_{p,m,h,j,t} = MSO_{p,m,h,j,t+TT2_j}, \forall p, m, h, j, t \tag{11}$$

$INH_{p,h,t}$ represents the sample $p$ that arrives at hospital $h$ at time $t$:

$$INH_{p,h,t} = \sum_{m,j} FTD_{p,m,h,j,t}, \forall p, h, t \tag{12}$$

*Capacity constraints.* Eq. (13) calculates the capacity ($CAP_{m,t}$) of each of the manufacturing facilities $m$ at every time $t$, while Eq. (14) ensures that the therapies in manufacturing do not exceed the available capacity:

$$CAP_{m,t} = FCAP_m - \sum_{p,\hat{t}} INM_{p,m,\hat{t}}, \forall p, m, t, t - TMFE \leq \hat{t} \leq t \tag{13}$$

$$\sum_p INM_{p,m,t} - \sum_p OUTM_{p,m,t} \leq CAP_{m,t}, \forall p, m, t \tag{14}$$

*Network structure constraints.* Eq. (15)-(16) ensure that matches are established only between existing manufacturing facilities. $E1_m$ is a binary variable that takes the value of 1 if a manufacturing facility is established, while it is equal to 0 otherwise. Variables $X1_{c,m}$ and $X2_{m,h}$ are binary variables that can take the value of 1 if only if $E1_m = 1$.

$$X1_{c,m} \leq E1_m, \forall c, m \tag{15}$$

$$X2_{m,h} \leq E1_m, \forall c, m \tag{16}$$

Equations (17)-(18) ensure that only one transport mode $j$ from every therapy $p$ at every journey can be selected. Variables $Y1_{p,c,m,j,t}$ and $Y2_{p,m,h,j,t}$ are binary variables that take the value of 1 if a transport mode $j$ is selected for the transportation of therapy $p$ between two facilities at time $t$.

$$\sum_{c,m,j,t} Y1_{p,c,m,j,t} = 1, \forall p, c, m, j, t \tag{17}$$

$$\sum_{m,h,j,t} Y2_{p,m,h,j,t} = 1, \forall p, m, h, j, t \tag{18}$$

Equation (19) adds an upper bound in the total number of manufacturing facilities that can be established.

$$\sum_m E1_m \leq U^M \tag{19}$$

*Demand satisfaction.* The total rate of flow of each therapy $p$ arriving at hospital $h$ must be equal to the corresponding demand (Eq. 20):

$$\sum_{p,h,t} INH_{p,h,t} = NP\ , \forall p, h, t \tag{20}$$

*Logical constraints for transportation flows.* Therapies $p$ can be transported from clinical site $c$ to a manufacturing site $m$ (Eq. 21) and from a manufacturing site $m$ to a hospital $h$ (Eq. 22) if and only if a match between the corresponding facilities has been previously established.

$$Y1_{p,c,m,j,t} \leq X1_{c,m}, \forall p, c, m, j, t \tag{21}$$

$$Y2_{p,m,h,j,t} \leq X2_{m,h}, \forall p, m, q, j, t \tag{22}$$

Eq. (23)-(26) make sure that a match is only made between a leukapheresis site $c$ and its corresponding co-located hospital $h$.

$$\sum_{m,j,t} Y2_{p,m,h1,j,t} \leq \sum_t INC_{p,c1,t} \cdot t\ , \forall p \tag{23}$$

$$\sum_{m,j,t} Y2_{p,m,h2,j,t} \leq \sum_t INC_{p,c2,t} \cdot t,\ \forall p \tag{24}$$

$$\sum_{m,j,t} Y2_{p,m,h3,j,t} \leq \sum_t INC_{p,c3,t} \cdot t,\ \ \forall p \tag{25}$$

$$\sum_{m,j,t} Y2_{p,m,h4,j,t} \leq \sum_t INC_{p,c4,t} \cdot t,\ \ \forall p \tag{26}$$

Eq. (27)-(30) ensure that a minimum and maximum flow of material exists for a transportation link to be established. The values for $FMIN$ and $FMAX$ can be established by the method presented by Tsiakis et al. (Tsiakis et al., 2001). Nonetheless, in this case where every therapy $p$ corresponds to a single patient (material), $FMIN$ and $FMAX$ are assumed to be equal to 0 and 1 respectively.

$$LSR_{p,c,m,j,t} \geq Y1_{p,c,m,j,t} \cdot FMIN, \forall p, c, m, j, t \tag{27}$$

$$LSR_{p,c,m,j,t} \leq Y1_{p,c,m,j,t} \cdot FMAX, \forall p, c, m, j, t \quad (28)$$

$$MSOQ_{p,m,hj,t} \geq Y2_{p,m,h,j,t} \cdot FMIN, \forall p, m, h, j, t \quad (29)$$

$$MSOQ_{p,m,hj,t} \leq Y2_{p,m,h,j,t} \cdot FMAX, \forall p, m, h, j, t \quad (30)$$

*Time constraints.* Eq. (31) calculates the manufacturing time $t$ of therapy $p$ in facility $m$.

$$DURM_{p,m,t} = \sum_{\hat{t}}(INM_{p,m,\hat{t}-1} - OUTM_{p,m,\hat{t}}) + OUTM_{p,m,t}, \forall\, p, m, t, \hat{t} \leq t \quad (31)$$

Eq. (32) establishes the time point when a patient checks into a leukapheresis site $c$, while Eq. (33) presents the time point when therapy $p$ is delivered to hospital $h$.

$$CTT_p = \sum_{h,t} INH_{p,h,t} \cdot t\,, \forall p \quad (32)$$

$$STT_p = \sum_{c,t} INC_{p,c,t} \cdot t\,, \forall p \quad (33)$$

Eq. (34) makes sure that the time point a patient checks into a leukapheresis site $c$ chronologically precedes the time point the corresponding therapy $p$ is delivered to hospital $h$.

$$STT_p \leq CTT_p, \forall p \quad (34)$$

Eq. (42) ensures that the turnaround time of therapy $p$ is less than or equal to a specified number of days $(ND)$.

$$TRT_p \leq ND \quad (35)$$

Eq. (36) calculates the total time (vein-to-vein) for a therapy $p$ $(TRT_p)$ from the time point that a patient checks into a leukapheresis center $(STT_p)$ until the time point that the therapy for this patient is delivered to the hospital location $(CTT_p)$:

$$TRT_p = CTT_p - STT_p, \forall p$$

(36)

Eq. (37) presents the average return time $(ATRT)$ of all the therapies $p$.

$$ATRT = \frac{\sum_p TRT_p}{NP} \quad (37)$$

Table A1. Model nomenclature for the Mixed Integer Linear Programming supply chain optimization problem formulation.

| | | |
|---|---|---|
| **SETS/DIMENSIONS** | | |
| $\boldsymbol{c}$ | leukapheresis sites | |
| $\boldsymbol{h}$ | hospitals | |
| $\boldsymbol{j}$ | transport modes | |
| $\boldsymbol{m}$ | manufacturing sites | |
| $\boldsymbol{p}$ | Patients/therapies | |
| $\boldsymbol{t}$ | time periods | |
| PARAMETERS | | |
| $\boldsymbol{CIM_m}$ | Capital investment for manufacturing facility $m$ | USD |
| $\boldsymbol{CFVM_m}$ | Fixed variable manufacturing cost for manufacturing facility $m$ (personnel, facilities and equipment maintenance costs) | USD |
| $\boldsymbol{CVM_p}$ | Variable cost of materials required for therapy $p$ | USD |
| $\boldsymbol{CQC}$ | Quality control cost when in house | USD $\text{therapy}^{-1}$ |
| $\boldsymbol{FCAP_m}$ | Total capacity of manufacturing site $m$ | therapies |
| $\boldsymbol{FMAX}$ | Maximum flow | - |
| $\boldsymbol{FMIN}$ | Minimum flow | - |
| $\boldsymbol{INC_{p,c,t}}$ | Patient $p$ arriving for leukapheresis site $c$ at time $t$ | patients |
| $\boldsymbol{TLS}$ | Duration of leukapheresis procedure | days |
| $\boldsymbol{TMFE}$ | Duration of the manufacturing process | days |
| $\boldsymbol{TQC}$ | Duration of Quality Control | days |
| $\boldsymbol{NT}$ | Number of time periods | - |
| $\boldsymbol{NP}$ | Number of patients/therapies | - |
| $\boldsymbol{TT1_j}$ | Transport time from leukapheresis site $c$ to manufacturing facility $m$ using transport mode $j$ | days |
| $\boldsymbol{TT2_j}$ | Transport time from manufacturing facility $m$ to hospital $h$ using transport mode $j$ | days |
| $\boldsymbol{U1_{c,m,j}}$ | Unit transport cost from leukapheresis site $c$ to manufacturing facility $m$ using transport mode $j$ | USD therapy$^{-1}$ day$^{-1}$ |
| $\boldsymbol{U2_{m,h,j}}$ | Unit transport cost from manufacturing facility $m$ to hospital $h$ using transport mode $j$ | USD therapy$^{-1}$ day$^{-1}$ |
| BINARY VARIABLES | | |
| $\boldsymbol{E1_{m}}$ | 1 if manufacturing facility $m$ is established | - |
| $\boldsymbol{X1_{c,m}}$ | 1 if a match between a leukapheresis site $c$ and manufacturing facility $m$ is established | - |
| $\boldsymbol{X2_{m,h}}$ | 1 if a match between manufacturing facility $m$ and a hospital $h$ is established | - |
| $\boldsymbol{Y1_{p,c,m,j,t}}$ | 1 a sample $p$ is transferred from a leukapheresis site c to a manufacturing facility $m$ via mode $j$ at time $t$ | - |
| $\boldsymbol{Y2_{p,m,h,j,t}}$ | 1 a sample $p$ is transferred from a manufacturing facility $m$ to a hospital h via mode $j$ at time $t$ | - |
| VARIABLES | | |
| $\boldsymbol{ATRT}$ | Average return time | days |
| $\boldsymbol{CAP_{m,t}}$ | Capacity of manufacturing facility $m$ at time $t$ | therapies day$^{-1}$ |
| $\boldsymbol{CTM_{p}}$ | Total manufacturing cost of therapy $p$ | USD therapy$^{-1}$ |
| $\boldsymbol{CTT_{p}}$ | Completion time of treatment for patient $p$ | - |
| $\boldsymbol{DURM_{p,m,t}}$ | 1 only for the time points $t$ at which a therapy $p$ is manufactured in facility $m$ | - |
| $\boldsymbol{FTD_{p,m,h,j,t}}$ | Final therapy that left from manufacturing facility $m$ arriving at hospital $h$ via mode j at time $t$ | - |
| $\boldsymbol{FTR_{p,m,h,j,t}}$ | Therapy $p$ leaving manufacturing facility $m$ to hospital $h$ via mode $j$ at time $t$ | - |
| $\boldsymbol{INH_{p,h,t}}$ | Therapy $p$ arriving at hospital $h$ at time $t$ | - |
| $\boldsymbol{INM_{p,m,t}}$ | Therapy $p$ arriving at manufacturing facility $m$ at time $t$ | - |
| $\boldsymbol{LSA_{p,c,m,j,t}}$ | Therapy $p$ that left leukapheresis site $c$ and arrived at manufacturing facility $m$ via mode $j$ at time $t$ | - |
| $\boldsymbol{LSR_{p,c,m,j,t}}$ | Therapy $p$ that left leukapheresis site $c$ to go to manufacturing facility $m$ via mode $j$ at time $t$ | - |
| $\boldsymbol{MSO_{p,m,h,j,t}}$ | Therapy $p$ leaving manufacturing facility $m$ to go to hospital $h$ via mode $j$ at time $t$ | - |
| $\boldsymbol{OUTC_{p,c,t}}$ | Therapy $p$ leaving leukapheresis site $c$ at time $t$ | - |
| $\boldsymbol{OUTM_{p,m,t}}$ | Therapy $p$ leaving manufacturing facility $m$ at time $t$ | - |
| $\boldsymbol{RATIO_{m,t}}$ | Percentage of utilization of manufacturing site $m$ at time $t$ | - |
| $\boldsymbol{STT_p}$ | Starting time of treatment for patient $p$ | - |
| $\boldsymbol{TOTCOST}$ | Total supply chain cost | USD |
| $\boldsymbol{TRT_p}$ | Total return time of therapy | days |
| $\boldsymbol{TTC_p}$ | Total transport cost per therapy $p$ | USD therapy$^{-1}$ |

Table A2. Candidate leukapheresis sites, manufacturing facilities, and hospitals as considered in this work.

| | Facility | Location | Capacity (parallel lines) |
|---|---|---|---|
| **Manufacturing sites** | $m_1$ | UK/Stevenage | 4 |
| | $m_2$ | EU/Berlin | 31 |
| | $m_3$ | EU/Belgium | 10 |
| | $m_4$ | USA/Pennsylvania | 4 |
| | $m_5$ | USA/Virginia | 31 |
| | $m_6$ | UK/Glasgow | 10 |
| **Leukapheresis sites** | $c_1$ | London | 8 patients/day |
| | $c_2$ | Glasgow | |
| | $c_3$ | Manchester | |
| | $c_4$ | Birmingham | |
| **Hospitals** | $h_1$ | London | N/A |
| | $h_2$ | Glasgow | N/A |
| | $h_3$ | Manchester | N/A |
| | $h_4$ | Birmingham | N/A |

**Appendix B**

Given a test data set $DataSet = \{InstanceSet_i,\ y_{m,i}\}_{i=1}^{N}$, where $y_{m,i}\ \in \{0,1\}^{|M|+1}$ is the true label of the $i^{th}$ instance (i.e. sample) and $\hat{y}_{m,i}$ is its predicted label, the following evaluation metrics are defined (Huang et al., 2015). The following evaluation metrics are transformed for multi-label classification and can be categorized as either sample-based or label-based. Higher values indicate better classifier performance, with the exception of Hamming loss, which should be 0 for a perfect classifier.

*Hamming loss* (Eq. 1) accounts for partial correctness and evaluates the fraction of labels that are misclassified (i.e., the symmetric difference between the set of true labels and the set of predicted labels) by applying exclusive disjunction.

$$Hamming\ loss = \frac{1}{N}\sum_{i=1}^{N}\frac{1}{|M|+1}\sum_{j=1}^{|M|+1}(y_{m,i} \oplus \hat{y}_{m,i}) \tag{1}$$

*Jaccard index* (Eq. 2) also known as label-level accuracy or intersection over union compares the proportion of correct labels to the total number of labels.

$$Jaccard\ Index = \frac{1}{N}\sum_{i=1}^{N}\frac{|y_{m,i} \wedge \hat{y}_{m,i}|}{|y_{m,i} \vee \hat{y}_{m,i}|} \tag{2}$$

*Sample-level accuracy* (Eq. 3) or exact match ratio accounts only for total correctness and evaluates the proportion of samples that have all their labels classified correctly.

$$Sample\ level\ Accuracy = \frac{1}{N}\sum_{i=1}^{N}(y_{m,i} \equiv \hat{y}_{m,i}) \tag{3}$$

*Precision* is a measure of the accuracy of the positive predictions and *Recall* is a measure of the ability of the model to capture all the positive instances. High precision indicates that when the classifier predicts a positive class, it is likely to be correct, while high recall indicates that the classifier is effective at identifying positive instances. *Sample-based Precision* (Eq. 4) and *sample-based Recall* (Eq. 5) are calculated for each sample/instance.

$$Precision = \frac{1}{N}\sum_{i=1}^{N}\frac{|y_{m,i} \wedge \hat{y}_{m,i}|}{|\hat{y}_{m,i}|} \tag{4}$$

$$Recall = \frac{1}{N}\sum_{i=1}^{N}\frac{|y_{m,i} \wedge \hat{y}_{m,i}|}{|y_{m,i}|} \tag{5}$$

Macro-averaged and micro-averaged metrics provide different perspectives on the performance of the classifier. Micro-averaged metrics are sample-based and give equal weight to each sample/instance, making them suitable for imbalanced datasets. On the other hand, macro-averaged metrics are label-based and treat each label equally, which is useful when all labels are equally important. The choice between macro and micro metrics depends on the specific problem at hand and the objective of the classification. Below are given the equations for macro-averaged precision (Eq. 6), macro-averaged recall (Eq. 7), micro-averaged precision (Eq. 8), and micro-averaged recall (Eq. 9).

$$macro\ Precision = \frac{1}{|M|+1}\sum_{j=1}^{|M|+1}\frac{|y_{m,j}\ \wedge\ \hat{y}_{m,j}|}{|\ \hat{y}_{m,j}\ |} \tag{6}$$

$$macro\ Recall = \frac{1}{|M|+1}\sum_{j=1}^{|M|+1}\frac{|y_{m,j}\ \wedge\ \hat{y}_{m,j}|}{|\ y_{m,i}\ |} \tag{7}$$

$$micro\ Precision = \frac{\sum_{j=1}^{|M|+1}\sum_{i=1}^{N}|y_{m,i}\ \wedge\ \hat{y}_{m,i}|}{\sum_{j=1}^{|M|+1}\sum_{i=1}^{N}|\ \hat{y}_{m,i}\ |} \tag{8}$$

$$micro\ Recall = \frac{\sum_{j=1}^{|M|+1}\sum_{i=1}^{N}|y_{m,i}\ \wedge\ \hat{y}_{m,i}|}{\sum_{j=1}^{|M|+1}\sum_{i=1}^{N}|\ y_{m,i}\ |} \tag{9}$$

*F1 score* is a metric that combines both precision and recall into a single value, providing a balanced assessment of a model's performance, especially in situations where there is an imbalance between positive and negative classes. *Macro-averaged* and *micro-averaged F1 score* are simply the harmonic means of equations (6)-(7) and (8)-(9) respectively. A perfect *F1 score* is 1.

Table B1. Evaluation metrics for the feed-forward neural network in the benchmark case study.

| **Class** | **Precision** | **Recall** | **F1-score** | weight |
|---|---|---|---|---|
| $m_1$ | 0.77 | 0.75 | 0.76 | 99 |
| $m_2$ | 0.85 | 0.89 | 0.87 | 100 |
| $m_3$ | 0.77 | 0.73 | 0.75 | 101 |
| $m_4$ | 0.95 | 0.86 | 0.90 | 100 |
| $m_5$ | 0.91 | 0.90 | 0.91 | 101 |
| $m_6$ | 0.93 | 0.79 | 0.85 | 99 |
| infeasible | 0.94 | 0.91 | 0.92 | 100 |
| micro avg | 0.83 | 0.84 | 0.83 | 700 |
| macro avg | 0.87 | 0.84 | 0.85 | 700 |
| weighted avg | 0.87 | 0.84 | 0.85 | 700 |

Table B2. Evaluation metrics for the convolutional neural network in the benchmark case study.

| **Class** | **Precision** | **Recall** | **F1-score** | weight |
|---|---|---|---|---|
| $m_1$ | 0.87 | 0.88 | 0.88 | 101 |
| $m_2$ | 0.94 | 0.96 | 0.95 | 100 |
| $m_3$ | 0.94 | 0.88 | 0.91 | 99 |
| $m_4$ | 0.98 | 0.94 | 0.96 | 100 |
| $m_5$ | 0.95 | 0.96 | 0.96 | 100 |
| $m_6$ | 0.97 | 0.90 | 0.93 | 101 |
| infeasible | 0.97 | 0.96 | 0.96 | 99 |
| micro avg | 0.93 | 0.93 | 0.93 | 700 |
| macro avg | 0.95 | 0.93 | 0.94 | 700 |
| weighted avg | 0.95 | 0.93 | 0.94 | 700 |

## Appendix C

Table C1. Feed-forward neural network predictions for the benchmark case study compared to the actual global optimum solution in the 12 instances from the validation set.

| **Demand** | **Distribution** | **Actual solution** | **Prediction** (50% probability threshold) | **Prediction** (1% probability threshold) |
|---|---|---|---|---|
| 200 | Uniform | $m_1, m_4$ | $m_1, m_4$ | $m_1, m_3, m_4$ |
| | Left triangular | $m_3$ | $m_3$ | $m_1, m_3, m_4$ |
| | Right triangular | $m_1, m_3$ | $m_1, m_4$ | $m_1, m_3, m_4$ |
| 500 | Uniform | $m_3, m_6$ | $m_3, m_6$ | $m_1, m_3, m_6$ |
| | Left triangular | $m_2$ | $m_2$ | $m_2, m_3, m_6$ |
| | Right triangular | $m_2$ | $m_3, m_6$ | $m_1, m_3, m_6$ |
| 1000 | Uniform | $m_2$ | $m_1, m_2$ | $m_2, m_3, m_6$ |
| | Left triangular | $m_2$ | $m_2$ | $m_1, m_2, m_3$ |
| | Right triangular | $m_2$ | $m_2$ | $m_1, m_2$ |
| 2000 | Uniform | $m_2, m_5$ | $m_2, m_5$ | $m_2, m_5$ |
| | Left triangular | infeasible | infeasible | infeasible |
| | Right triangular | infeasible | infeasible | infeasible |

Table C2. Convolutional neural network predictions for the benchmark case study compared to the actual global optimum solution in the 12 instances from the validation set.

| **Demand** | **Distribution** | **Actual solution** | **Prediction** (50% probability threshold) | **Prediction** (1% probability threshold) |
|---|---|---|---|---|
| 200 | Uniform | $m_1, m_4$ | $m_1, m_4$ | $m_1, m_3, m_4$ |
| | Left triangular | $m_3$ | $m_3$ | $m_1, m_3$ |
| | Right triangular | $m_1, m_3$ | $m_3$ | $m_1, m_3, m_4$ |
| 500 | Uniform | $m_3, m_6$ | $m_3, m_6$ | $m_2, m_3, m_6$ |
| | Left triangular | $m_2$ | $m_2$ | $m_2$ |
| | Right triangular | $m_2$ | $m_1, m_3$ | $m_1, m_2, m_3, m_6$ |
| 1000 | Uniform | $m_2$ | $m_2$ | $m_2$ |
| | Left triangular | $m_2$ | $m_2$ | $m_1, m_2$ |
| | Right triangular | $m_2$ | $m_1, m_2$ | $m_1, m_2$ |
| 2000 | Uniform | $m_2, m_5$ | $m_2, m_5$ | $m_2, m_5$ |
| | Left triangular | infeasible | infeasible | infeasible |
| | Right triangular | infeasible | $m_2, m_5$ | $m_2, m_5$ |